\documentclass[12pt,a4paper]{amsart}
\usepackage{graphicx}
\usepackage{amssymb}
\usepackage{amsfonts}
\usepackage{epsfig}
\usepackage{amscd}
\usepackage{pstricks}

\newlength{\cellsz}
\newcounter{cellsize}
\newcommand{\setcellsize}[1]{%
  \setcounter{cellsize}{#1}%
  \setlength{\cellsz}{\value{cellsize}\unitlength}}%
\setcellsize{13}%
\newcommand\cellify[1]{\def\thearg{#1}\def\nothing{}%
\ifx\thearg\nothing
\vrule width0pt height\cellsz depth0pt\else
\hbox to 0pt{{\begin{picture}(\value{cellsize},\value{cellsize})
  \put(0,0){\line(1,0){\value{cellsize}}}
  \put(0,0){\line(0,1){\value{cellsize}}}
  \put(\value{cellsize},0){\line(0,1){\value{cellsize}}}
  \put(0,\value{cellsize}){\line(1,0){\value{cellsize}}} \end{picture} \hss}}\fi%
\vbox to \cellsz{ \vss \hbox to \cellsz{\hss$#1$\hss} \vss}}
\newcommand\tableaux[1]{\vcenter{\vbox{\let\\\cr
\baselineskip -16000pt \lineskiplimit 16000pt \lineskip 0pt
\ialign{&\cellify{##}\cr#1\crcr}}}}
\newcommand\tabl[1]{\vtop{\let\\\cr
\baselineskip -16000pt \lineskiplimit 16000pt \lineskip 0pt
\ialign{&\cellify{##}\cr#1\crcr}}}

\edef\savecatcodeat{\the\catcode`@}
\catcode`\@=11

\def\tb@ifSpecChars#1#2{#1}
\def\tb@ifNoSpecChars#1#2{#2}

\def\tableau{%
  \bgroup% matched in \tb@tableauD
  \@ifstar{\let\Tif\tb@ifNoSpecChars\tb@tableauB}% *, don't use special chars
          {\let\Tif\tb@ifSpecChars\tb@tableauB}}% no *, use special chars

\def\tb@tableauB{% add [] if no [options]
  \@ifnextchar[{\tb@tableauC}{\tb@tableauC[]}}

\def\tb@tableauC[#1]{\hbox\bgroup%
    \let\\=\cr% end line
    \def\bl{\global\let\tbcellF\tb@cellNF}%
    \def\tf{\global\let\tbcellF\tb@cellH}% highlighted cell
    \dimen2=\ht\strutbox \advance\dimen2 by\dp\strutbox%
    \ifx\baselinestretch\undefined\relax%
    \else%
       \dimen0=100sp \dimen0=\baselinestretch\dimen0%
       \dimen2=100\dimen2 \divide\dimen2 by\dimen0%
    \fi%
    \let\tpos\tb@vcenter% default position
    \tb@initYoung% default tableau type
    \tb@options#1\eoo% parse options
    \let\arrow\tb@arrow%
    \dimen0=\Tscale\dimen2%
    \dimen1=\dimen0 \advance\dimen1 by \tb@fframe%
    \lineskip=0pt\baselineskip=0pt% line spacing will be from \vbox to \dimen0
    \def\tb@nothing{}%
    \def\endcellno{$\rss\egroup\bss\egroup}% end cell w/o overlap
    \def\endcell{\endcellno\kern-\dimen0}% end cell & prepare to overlap it
    \def\begincell{\vbox to\dimen0\bgroup\vss\hbox to\dimen0\bgroup\hss$}%
    \let\overlay\tb@overlay%
    \let\fl\tb@fl%
    \let\lss\hss\let\rss\hss\let\tss\vss\let\bss\vss% cell alignment
    \def\mkcell##1{% format individual cell
        \let\tbcellF\tb@cellD% default cell frame
        \def\tb@cellarg{##1}% store cell contents
        \ifx\tb@cellarg\tb@nothing\let\tb@cellarg\tb@cellE\fi%
        \begincell\tb@cellarg\endcellno% the actual cell content
        \tbcellF}% draw cell frame
    \let\savecellF\tbcellF% save global value of cellF in case of nested tableau
     \Tif{\catcode`,=4\catcode`|=\active}{}\tb@tableauD}%

\let\tb@savetableauD\tableauD% save any current definition
{% set up characters which will be interpreted as command characters
    \catcode`|=\active \catcode`*=\active \catcode`~=\active%
    \catcode`@=\active% command characters
\gdef\tableauD#1{%
  \Tif{% make all the command characters active in math mode when #1 parsed
    \mathcode`|="8000 \mathcode`*="8000%
    \mathcode`~="8000 \mathcode`@="8000%
    \def@{\bullet}%
    \let|\cr% end line
    \let*\tf% highlighted cell
    \let~\sk% skew cell
  }{}%
  \tpos{\tabskip=0pt\halign{&\mkcell{##}\cr#1\crcr}}%
  \global\let\tbcellF\savecellF% restore global value
  \egroup% match \hbox\bgroup at start of \tableauC
  \egroup}% match \bgroup at start of \tableau
}
\let\tb@tableauD\tableauD% rename the command
\let\tableauD\tb@savetableauD% restore old command with this name
\let\tb@savetableauD\undefined

\def\tb@options#1{\ifx#1\eoo\relax\else\tb@option#1\expandafter\tb@options\fi}

\def\tb@option#1{%
  \if#1t\let\tpos\tb@vtop\fi%        t = align at top
  \if#1c\let\tpos\tb@vcenter\fi%     c = align at center
  \if#1b\let\tpos\vbox\fi%           b = align at bottom
  \if#1F\tb@initFerrers\fi%          F = Ferrers diagram
  \if#1Y\tb@initYoung\fi%            Y = Young diagram
  \if#1s\tb@initSmall\fi%            s = small boxes
  \if#1m\tb@initMedium\fi%           m = medium boxes
  \if#1l\tb@initLarge\fi%            l = large boxes
  \if#1p\tb@initPartition\fi%            p = small partition sized boxes
  \if#1a\tb@initArrow\fi%            a = use arrow font as base dimension
}

\def\tb@vcenter#1{\ifmmode\vcenter{#1}\else$\vcenter{#1}$\fi}

\def\tb@vtop#1{\hbox{\raise\ht\strutbox\hbox{\lower\dimen0\vtop{#1}}}}

\def\tb@initPartition{\def\Tscale{.3}}
\def\tb@initSmall{\def\Tscale{1}}
\def\tb@initMedium{\def\Tscale{2}}
\def\tb@initLarge{\def\Tscale{3}}

\def\tb@initArrow{\dimen2=1.25em}

\def\tb@initYoung{%
  \def\tb@cellE{}% empty cell stays empty
  \let\tb@cellD\tb@cellN% default frame is normal frame
  \def\sk{\global\let\tbcellF\tb@cellNF}}% skew cells are empty
\def\tb@initFerrers{%
  \def\tb@cellE{\bullet}% empty cell gets bullet
  \let\tb@cellD\tb@cellNF% default frame is no frame
  \def\sk{\bullet}}% skew cell gets bullet

\tb@initMedium% default scale

\def\tb@sframe#1{%
  \vbox to0pt{%            Embed frame in a box of no vert or hor extent
    \vss%                            pull box above reference point
    \hbox to0pt{%
      \hss%                          pull box left of reference point
      \vbox to\dimen1{%              Actual width of frame
        \hrule depth #1 height0pt% draw top edge of frame
        \vss%                     begin vcenter sides
        \hbox to\dimen1{%           horiz box with side edges just inside
          \vrule width #1 height\dimen1% left edge
          \hss%                     stretch center
          \vrule width #1%         right edge
          }%
        \vss%                     end vcenter sides
        \hrule height #1 depth 0in% bottom edge
        }%
      \kern-\tb@hframe%           horiz alignment off by half line width
      }%
    \kern-\tb@hframe}}%           vert alignment off by half line width
\def\tb@hframe{.2pt}\def\tb@fframe{.4pt}\def\tb@bframe{2pt}
\def\tb@cellH{\tb@sframe{\tb@bframe}}       % bold frame
\def\tb@cellNF{}                            % no frame
\def\tb@cellN{\tb@sframe{\tb@fframe}}       % normal frame
\let\tbcellF\tb@cellN                       % default is normal

\def\tb@rpad{1pt}
\def\tb@lpad{1pt}
\def\tb@tpad{1.8pt}
\def\tb@bpad{1.8pt}

\def\tb@overlay{\endcell\@ifnextchar[{\tb@overlaya}{\begincell}}
\def\tb@overlaya[#1]{\vbox to\dimen0\bgroup%
  \tb@overlayoptions#1\eoo%
  \tss\hbox to\dimen0\bgroup\lss$}
\def\tb@overlayoptions#1{\ifx#1\eoo\relax\else\tb@overlayoption#1\expandafter\tb@overlayoptions\fi}

\def\tb@overlayoption#1{
  \if#1t\def\tss{\vskip\tb@tpad}\let\bss\vss\fi% t = align at top
  \if#1c\let\tss\vss\let\bss\vss\fi%             c = align at center
  \if#1b\def\bss{\vskip\tb@bpad}\let\tss\vss\fi% b = align at bottom
  \if#1l\def\lss{\hskip\tb@lpad}\let\rss\hss\fi% l = align at left
  \if#1m\let\lss\hss\let\rss\hss\fi%             m = align at middle
  \if#1r\def\rss{\hskip\tb@rpad}\let\lss\hss\fi% r = align at right
}

\def\tb@fl{\endcell\begincell\vrule depth 0pt width \dimen0 height \dimen0 \endcell\begincell}

\@ifundefined{diagram}{}{
\def\tb@arrowpad{.5}

\newoptcommand{\tb@arrow}{\@ne}[2]{%
  \endcell% end previous cell contents
   \begingroup%
   \let\dg@getnodesize\tb@getnodesize% substitute routine to get nodesize
   \dg@USERSIZE=#1\relax%
   \ifnum\dg@USERSIZE<\@ne \dg@USERSIZE=\@ne \fi%
   \dg@parse{#2}%
   \dg@label{\tb@draw{#1}{#2}}}% draw arrow

\def\tb@getnodesize#1#2#3#4#5{\dimen3=\tb@arrowpad\dimen2 #4=\dimen3 #5=\dimen3\relax}
\def\tb@getnodesize#1#2#3#4#5{\ifnum#2=0\ifnum#3=0\tb@getnodesizetail{#4}{#5}\else\tb@getnodesizehead{#4}{#5}\fi\else\tb@getnodesizehead{#4}{#5}\fi}
\def\tb@getnodesizetail#1#2{\dimen3=.5\dimen2 #1=\dimen3 #2=\dimen3}
\def\tb@getnodesizehead#1#2{\dimen3=.5\dimen2 #1=\dimen3 #2=\dimen3}

\def\tb@draw#1#2#3#4{%
        \dg@X=0\dg@Y=0\dg@XGRID=1\dg@YGRID=1\unitlength=.001\dimen0%
        \dg@LBLOFF=\dgLABELOFFSET \divide\dg@LBLOFF\unitlength%
        \dg@drawcalc% compute arrow geometry
        \begincell% start tableau cell
        \let\lams@arrow\tb@lams@arrow% substitute routine
        \begin{picture}(0,0)\begingroup\dg@draw{#1}{#2}{#3}{#4}\end{picture}%
        \endcell% end tableau cell
        \endgroup% match \begingroup in \tb@arrow
        \begincell}% start new entry in this cell
}

\def\tb@lams@arrow#1#2{%
 \lams@firstx\z@\lams@firsty\z@
 \lams@lastx#1\relax\lams@lasty#2\relax
 \lams@center\z@
 \N@false\E@false\H@false\V@false
 \ifdim\lams@lastx>\z@\E@true\fi
 \ifdim\lams@lastx=\z@\V@true\fi
 \ifdim\lams@lasty>\z@\N@true\fi
 \ifdim\lams@lasty=\z@\H@true\fi
 \NESW@false
 \ifN@\ifE@\NESW@true\fi\else\ifE@\else\NESW@true\fi\fi
 \ifH@\else\ifV@\else
  \lams@slope
  \ifnum\lams@tani>\lams@tanii
   \lams@ht\ten@\p@\lams@wd\ten@\p@
   \multiply\lams@wd\lams@tanii\divide\lams@wd\lams@tani
  \else
   \lams@wd\ten@\p@\lams@ht\ten@\p@
   \divide\lams@ht\lams@tanii\multiply\lams@ht\lams@tani
  \fi
 \fi\fi
 \ifH@  \lams@harrow
 \else\ifV@ \lams@varrow
 \else \lams@darrow
 \fi\fi
}

\catcode`\@=\savecatcodeat
\let\savecatcodeat\undefined

\newif\ifcoprod
\coprodfalse

\ifcoprod
\fi

\def\PT{{\rm PT}}

\def\Lh{{\rm Lh}}
\def\PP{{\rm PP}}
\def\tto{{\rm tot}}
\def\ttog{{\rm totg}}
\def\stb{{\rm st'}}

\def\Le{\hbox{\rotatedown{$\Gamma$}}}
\def\st{{\rm st}}
\def\tS{\tilde S}
\def\maj{{\rm maj}}
\def\mytr{{\,\triangleright}}
\def\Lig{{\rm LC}}
\def\MM{{\mathcal M}}
\def\VV{{\mathcal V}}
\def\LL{{\mathcal F}}
\def\ss{\scriptstyle}
\def\empd#1#2{{\scriptstyle \begin{matrix}\ss #1 \cr \ss #2\end{matrix}}}

\newtheorem{example}{Example}[section]

\newtheorem{theorem}[example]{Theorem}

\newtheorem{corollary}[example]{Corollary}
\newtheorem{conjecture}[example]{Conjecture}

\newtheorem{proposition}[example]{Proposition}

\newtheorem{lemma}[example]{Lemma}

\def\Proof{\noindent \it Proof -- \rm}
\def\qed{\hspace{3.5mm} \hfill \vbox{\hrule height 3pt depth 2 pt width 2mm}
\bigskip}

\def\sinv{{\rm sinv}}
\def\MM{{\mathcal M}}  % Base standard de TD.
\def\TEb{{\bf T'}}

\def\M{{\bf M}}       % Base standard de WQSym.

\def\GDes{\operatorname{GDes}}
\def\GC{\operatorname{GC}}
\def\WC{\operatorname{WC}}
\def\Col{{\rm Col}}
\def\raff{\succeq}     % Ordre de raffinement $I\raff J$ si I + fin que J.
\def\raffi{\preceq}    % Ordre de raffinement $I\raffi J$ si I - fin que J.
\def\finer{\raff}      % Ordre de raffinement $I\raffi J$ si I - fin que J.
\def\fatter{\raffi}    % Ordre de raffinement $I\raffi J$ si I - fin que J.

\def\limproj{{\rm proj\,lim}}
\def\K{{\mathbb K}}   % corps de base
\def\Sym{{\bf Sym}}   % NCSF
\def\qp{\star_q}
\def\NCSF{{\bf Sym}}           % NCSF
\def\FQSym{{\bf FQSym}}        % permutations
\def\MQSym{{\bf MQSym}}        % matrices
\def\WQSym{{\bf WQSym}}        % Mots initiaux
\def\C{{\mathbb C}}
\def\pack{{\rm pack}}          % packed word
\def\ev{{\rm ev}}       % evaluation
\def\inv{{\rm inv}}     % inversions
\def\ssh{\Cup}          % shuffle shifte
\def\sconc{\bullet}     % concatenation shiftee
\def\Std{{\rm Std}}     % standardisation
\def\convW{{*_W}}       % une loi de convolution
\def\<{\langle}
\def\>{\rangle}
\def\NN{{\mathbb N}}    % entiers naturels
\def\F{{\bf F}}         % F de FQSym, PQSym, ...
\def\G{{\bf G}}         % G de FQSym^*
\def\SG{{\mathfrak S}}  % groupe symetrique

\def\tensor{\otimes}

\def\Des{\operatorname{Des}}
\def\DesC{\operatorname{DC}}

\def\PW{{\rm PW}}   % packed words = mots tasses

\def\qbin#1#2{{\begin{bmatrix} #1 \\ #2\end{bmatrix}}_q}

\def\shuff#1#2{\mathbin{
      \hbox{\vbox{
        \hbox{\vrule
              \hskip#2
              \vrule height#1 width 0pt
               }%
        \hrule}%
             \vbox{
        \hbox{\vrule
              \hskip#2
              \vrule height#1 width 0pt
               \vrule }%
        \hrule}%
}}}
\def\shuffl{{\mathchoice{\shuff{7pt}{3.5pt}}%
                        {\shuff{6pt}{3pt}}%
                        {\shuff{4pt}{2pt}}%
                        {\shuff{3pt}{1.5pt}}}}%
\def\shuffle{\, \shuffl \,}

\newcommand{\T}{\mathcal{T}}
\newcommand{\pP}{\operatorname{\cdot}}

\newcommand{\mm}{\operatorname{\bf m}}

\newcommand{\Cr}{{cr}}

\def\rk{\mathrm{rank}}
\def\Stat{\mathrm{Fact}}
\def\QStat{\mathrm{QFact}}

\def\ash{\, \underline{\shuffl} \, }

\def\DessinMatrix#1{\vcenter{\hbox{\makebox[1.7ex]{$#1$}}}}
\def\GenMatrix#1{\vcenter{\halign{&$\DessinMatrix{##}$\cr#1}}\egroup}

\def\setinterlineskip#1{\baselineskip=0pt
  \lineskip=#1 \lineskiplimit=\maxdimen}

\def\matrice{%
  \bgroup
  \let\ =\omit
  \let\\=\cr
  \setinterlineskip{4.0pt}
  \GenMatrix}

\def\DessinsMatrix#1{\vcenter{\hbox{\makebox[1.3ex]{$\scriptstyle#1$}}}}
\def\GensMatrix#1{\vcenter{\halign{&$\DessinsMatrix{##}$\cr#1}}\egroup}

\def\MS{\mathbf{MS}}

\def\smallmatrice{%
  \bgroup
  \let\ =\omit
  \let\\=\cr
  \setinterlineskip{3.0pt}
  \GensMatrix}

\def\tensor{\otimes}

\newcommand{\SMat}[1]{{\left[\smallmatrice{#1\\}\right]}}

\newlength{\Hackl}
\newcommand{\Hack}{\vrule height \Hackl width 0pt}
\newcommand{\indexmat}%
    {\smallmatrice{\Hack a\\\Hack b\\\Hack c\\\Hack d\\\Hack e\\}}
\def\gf#1#2{\genfrac{}{}{0pt}{}{#1}{#2}}

\title[Hopf algebras and permutation tableaux]{Combinatorial Hopf algebras,
noncommutative Hall-Littlewood functions, and permutation tableaux}

\author[J.-C.~Novelli, J.-Y.~Thibon, and L. K. Williams]
{Jean-Christophe Novelli, Jean-Yves Thibon, and Lauren K. Williams}

\address[] {Institut Gaspard Monge, Universit\'e Paris-Est, Marne-la-Vall\'ee \\
5 Boulevard Des\-cartes \\Champs-sur-Marne \\77454 Marne-la-Vall\'ee cedex 2 \\
FRANCE}
\address[] {Harvard University, Department of Mathematics,
432 Science Center, 1 Oxford Street, Cambridge, MA 02138, USA}
\email[Jean-Christophe Novelli]{novelli@univ-mlv.fr}
\email[Jean-Yves Thibon]{jyt@univ-mlv.fr} 
\email[Lauren K. Williams]{lauren@math.harvard.edu}
\date{}

\begin{document}

\begin{abstract}
We introduce a new family of noncommutative analogues of the Hall-Littlewood
symmetric functions. Our construction relies upon Tevlin's bases and
simple $q$-deformations of the classical combinatorial Hopf algebras.
We connect our new Hall-Littlewood functions to permutation tableaux, and 
also give an exact formula for the $q$-enumeration of permutation 
tableaux of a fixed shape.  This gives an explicit formula for: 
the steady state probability of each state in the partially asymmetric
exclusion process (PASEP); the polynomial enumerating permutations with a
fixed set of {\it weak excedances} according to {\it crossings};
the polynomial enumerating permutations with a fixed set of 
{\it descent bottoms} according to occurrences of the {\it generalized
pattern} $2-31$.
\end{abstract}

\maketitle
\setcounter{tocdepth}{1}
\tableofcontents

%%%%%%%%%%%%%%%%%%%%%%%%%%%%%%%%%%%%%%%%%%%%%%%%%%%%%%%%%%%%%%%%%%%%%%%%%%%%%%%
%%%%%%%%%%%%%%%%%%%%%%%%%%%%%%%%%%%%%%%%%%%%%%%%%%%%%%%%%%%%%%%%%%%%%%%%%%%%%%%
%%%%%%%%%%%%%%%%%%%%%%%%%%%%%%%%%%%%%%%%%%%%%%%%%%%%%%%%%%%%%%%%%%%%%%%%%%%%%%%
\section{Introduction}

The combinatorics of Hall-Littlewood functions is one of the most interesting 
aspects of the modern theory of symmetric functions~\cite{Mcd}.
These are bases of symmetric functions, depending on a parameter $q$ which was
originally regarded as the cardinality of a finite field. They are named after
Littlewood's explicit realization of the Hall algebra in terms of symmetric
functions, which gave meaning to arbitrary complex values of
$q$~\cite{Li1}.

Combinatorics entered the scene with the observation by Foulkes~\cite{Fo1}
that the transition matrices between Schur functions and Hall-Littlewood
$P$-functions seemed to be given by polynomials, which were nonnegative
$q$-analogues of the well-known Kostka numbers, counting Young tableaux
according to shape and weight.
The conjecture of Foulkes was established by Lascoux and
Sch\"ut\-zen\-berger~\cite{LS2}, who introduced the charge statistic on Young
tableaux to explain the powers of $q$. Almost simultaneously,
Lusztig~\cite{Lu2} obtained an interpretation in terms of the intersection
homology of nilpotent orbits, and it is now known that these Kostka-Foulkes
polynomials are particular Kazhdan-Lusztig polynomials associated with the
affine Weyl groups of type $A$~\cite{Lus}.
But this was not the end of the story. Some ten years later, Kirillov and
Reshetikhin~\cite{KR} discovered an interpretation of Kostka-Foulkes
polynomials in statistical physics, as generating functions of Bethe ansatz
configurations for some generalizations of Heisenberg's $XXX$-magnet model,
and obtained a closed expression in the form of a sum of products of
$q$-binomial coefficients.

All of these results have been generalized in many directions. Generalized
Hall algebras (associated with quivers) have been introduced~\cite{Ring}.
Ribbon tableaux~\cite{LLT2} and $k$-Schur functions~\cite{LLM} give rise to
generalizations of the charge polynomials, sometimes interpretable as
Kazhdan-Lusztig polynomials~\cite{LT}.
Intersection homology has been computed for other varieties.
The Kirillov-Reshetikhin formula is now included in a vast corpus of fermionic
formulas, available for a large number of models~\cite{HKOTZ}.
However, the relations -if any- between these theories are generally unknown.

The present article is devoted to a different kind of generalization of
the Hall-Littlewood theory. It is by now well-known that many aspects of
the theory of symmetric functions can be lifted to Noncommutative symmetric
functions, or quasi-symmetric functions, and that those points which do
not have a good analogue at this level can sometimes be explained  by lifting
them to more complicated combinatorial Hopf algebras\footnote{
There is no general agreement on the precise definition of a combinatorial
Hopf algebra, see \cite{ABS} and \cite{LR-cha} for attempts at making this
concept precise, and \cite{NT, NT07, NTqthooks} for more examples.}.
The paradigm here is the
Littlewood-Richardson rule, which becomes trivial in the algebra of
Free symmetric functions, all the difficulty having 
been diluted in the definition of the algebra~\cite{NCSF6}.

A theory of noncommutative and quasi-symmetric Hall-Littlewood 
functions has been worked out
by Hivert~\cite{Hiv-adv}. Since there is no Hall algebra to use as a starting
point, Hivert's choice was to imitate Littlewood's definition, which can be
reformulated in terms of an action of the affine Hecke algebra on polynomials.
By replacing the usual action by a quasi-symmetrizing one, Hivert obtained
interesting bases, behaving in much the same way as the original ones, and
were easily deformable with a second parameter, so that analogues of
Macdonald's functions could also be defined \cite{HLT}.
See also the interesting work of Bergeron and Zabrocki \cite{BZ}.

However, Hivert's analogues of the Kostka-Foulkes polynomials are just
Kostka-Foulkes {\em monomials}, i.e., powers of $q$, given moreover by
a simple explicit formula. So the combinatorial connections to
tableaux, geometry and statistical physics do not show up in this
theory.

More recently, new possibilities arose with Tevlin's~\cite{Tev} discovery of
a plausible analogue of monomial symmetric functions on the noncommutative
side.
Tevlin's constructions are incompatible with the Hopf structure
(his monomial functions are not dual to products of complete functions
in any reasonable sense), so it seemed unlikely that they could lead to
interesting combinatorics. Nevertheless, Tevlin computed analogues of the
Kostka matrices in his setting, and conjectured that they had nonnegative
integer coefficients. This conjecture was proved in~\cite{HNTT},
and turned out to be more interesting than expected.
The proof required the use of larger combinatorial Hopf algebras, and
led to a vast generalization of the Gennocchi numbers.

In this paper we give a new generalization of Hall-Littlewood functions,
starting from Tevlin's bases. We define $q$-analogues $S^I(q)$ of the products
of complete homogeneous functions by embedding $\Sym$ in an associative
deformation of $\WQSym$ and projecting back to $\Sym$ by the map
introduced in \cite{HNTT}. This defines a nonassociative $q$-product
$\qp$ on $\Sym$, and our Hall-Littlewood functions are equal to the
products (see Section~\ref{sec-SL})
\begin{equation}
S^I(q) = S^{i_1} \qp (S^{i_2}\qp(\dots (S^{i_{r-1}}\qp S^{i_r}))).
\end{equation}

\noindent
These functions can be regarded as interpolating between the $S^I$ (at $q=1$)
and a new kind of noncommutative Schur functions (at $q=0$), have nonnegative
coefficients, which can be expressed in closed form as products of
$q$-binomial coefficients, and have a transparent combinatorial
interpretation.
As a consequence, the basis $R_I(q)$, defined by Moebius inversion on the
composition lattice, is also nonnegative on the same basis.

The really interesting phenomenon occurs with Tevlin's second basis
(denoted here and in~\cite{HNTT} by $L_I$),
an analogue of Gessel's fundamental basis $F_I$.
One can observe that the last column of the matrix $M(S,L)$ (which expresses
$S^{1^n}$ in terms of the $L_I$'s) gives the enumeration of
\emph{permutation tableaux} \cite{SW} by shape.
This observation can be easily shown, and one may wonder whether the expansion
of $S^{1^n}(q)$ on $L_I$ gives rise to interesting $q$-analogues. This is
clearly not the case (there are negative coefficients), but it turns out that
introducing a simple $q$-analogue $L_I(q)$ of $L_I$, we obtain again
nonnegative polynomials in the matrix $M(S(q),L(q))$. Finally, the matrix
$M(R(q),L(q))$ gives the $q$-enumeration of permutation tableaux according to
shape and rank, and another (yet unknown) statistic. Other (conjectural)
combinatorial interpretations in terms of permutations or packed words are
also proposed\footnote{Since many of our results and conjectures are stated in
terms of permutations, one might be tempted to work only with $\FQSym$,
bypassing $\WQSym$ entirely. However, as was already clear in~\cite{HNTT}, one
cannot make sense of the definition of Tevlin's monomial basis $\Psi$
using $\FQSym$.}.

As permutation tableaux occur in geometry (they are a distinguished subset
of Postnikov's {\em $\Le$-diagrams}, which parameterize cells in the totally
non-negative part of the Grassmannian \cite{Postnikov}) and the 
$q$-enumeration of permutation tableaux is (up to a shift) counting
cells according to dimension.  Additionally, permutation tableaux occur
in physics -- Corteel and Williams~\cite{CW} found a close connection to a 
well-known model from statistical physics called the 
partially asymmetric exclusion process, which in turn is related to the Hamiltonian 
of the  XXZ quantum spin chain \cite{ER}.
Therefore we may say that our new Hall-Littlewood
functions have some of the features which were absent from Hivert's theory.
However, we do not have the algebraic side coming from affine Hecke algebras,
and it is an open question whether both points of view can be unified.

We conclude this paper with exact formulas for the $q$-enumeration 
of permutation tableaux of types A and B, according to shape.
In the type A case, by the result of Corteel and Williams \cite{CW},
this gives an exact formula for the steady state probability of 
each state of the partially asymmetric exclusion process
(with arbitrary $q$ and $\alpha=\beta=\gamma=\delta=1$).  Applying
results of \cite{SW}, this also gives an exact formula for the number
of permutations with a fixed {\it weak excedance set} enumerated according to 
{\it crossings}, and for the number of permutations with a fixed set of 
{\it descent bottoms}, enumerated according to occurrences of the 
pattern $2-31$.

{\footnotesize {\it Acknowledgments.-}
This work has been partially supported by Agence Nationale de la Recherche,
grant ANR-06-BLAN-0380.
The authors would also like to thank the contributors of the MuPAD project,
and especially those of the combinat package, for providing the development
environment for this research (see~\cite{HT} for an introduction to
MuPAD-Combinat).
}

%%%%%%%%%%%%%%%%%%%%%%%%%%%%%%%%%%%%%%%%%%%%%%%%%%%%%%%%%%%%%%%%%%%%%%%%%%%%%%%
%%%%%%%%%%%%%%%%%%%%%%%%%%%%%%%%%%%%%%%%%%%%%%%%%%%%%%%%%%%%%%%%%%%%%%%%%%%%%%%
%%%%%%%%%%%%%%%%%%%%%%%%%%%%%%%%%%%%%%%%%%%%%%%%%%%%%%%%%%%%%%%%%%%%%%%%%%%%%%%
\section{Notations and background}

%%%%%%%%%%%%%%%%%%%%%%%%%%%%%%%%%%%%%%%%%%%%%%%%%%%%%%%%%%%%%%%%%%%%%%%%%%%%%%%
\subsection{Words, permutations, and compositions}

We assume that the reader is familiar with the standard notations of the
theory of noncommutative symmetric functions~\cite{NCSF1,NCSF6}.
%%%
We shall need an infinite totally ordered alphabet
$A=\{a_1<a_2<\cdots<a_n<\cdots\}$, generally assumed to be the set of
positive integers.
We denote by $\K$ a field of characteristic $0$, and by $\K\<A\>$ the free
associative algebra over $A$ when $A$ is finite, and the projective limit
$\limproj_B \K\<B\>$, where $B$ runs over finite subsets of $A$, when $A$ is
infinite.
The \emph{evaluation} $\ev(w)$ of a word $w$ is the sequence whose $i$-th term
is the number of times the letter $a_i$ occurs in $w$.
The \emph{standardized word} $\Std(w)$ of a word $w\in A^*$ is the permutation
obtained by iteratively scanning $w$ from left to right, and labelling
$1,2,\ldots$ the occurrences of its smallest letter, then numbering the
occurrences of the next one, and so on.
For example, $\Std(bbacab)=341625$.
For a word $w$ on the alphabet $\{1,2,\ldots\}$, we denote by $w[k]$ the word
obtained by replacing each letter $i$ by the integer $i+k$.
If $u$ and $v$ are two words, with $u$ of length $k$, one defines
the \emph{shifted concatenation}
$u\sconc v = u\cdot (v[k])$
and the \emph{shifted shuffle}
$ u\ssh v= u\shuffle (v[k])$,
where $\shuffle$ is the usual shuffle product.

Recall that a permutation $\sigma$ admits a {\it descent} at position
$i$ if $\sigma(i)> \sigma(i+1)$.  Symmetrically, $\sigma$
admits a {\it recoil} at $i$ if $\sigma^{-1}(i)>\sigma^{-1}(i+1)$.
The descent and recoil sets of $\sigma$ are the positions
of the descents and recoils, respectively.

A \emph{composition} of an integer $n$ is a sequence
$I=(i_1,\dots,i_r)$ of positive integers of sum $n$. 
In this case we write $I \models n$.  The integer $r$ is called
the \emph{length} of the composition.
The \emph{descent set} of $I$ is
$\Des(I) = \{ i_1,\ i_1+i_2, \ldots , i_1+\dots+i_{r-1}\}$.
The \emph{reverse refinement order}, denoted by $\raff$, on
compositions is such that $I=(i_1,\ldots,i_k)\raff J=(j_1,\ldots,j_l)$ iff
$\Des(I)\supseteq\Des(J)$, or equivalently,
$\{i_1,i_1+i_2,\ldots,i_1+\cdots+i_k\}$ contains
$\{j_1,j_1+j_2,\ldots,j_1+\cdots+j_l\}$.
In this case, we say that $I$ is finer than $J$.
For example, $(2,1,2,3,1,2)\raff (3,2,6)$.
The \emph{descent composition} $\DesC(\sigma)$ of a permutation
$\sigma\in\SG_n$ is the composition $I$ of $n$ whose descent set is the
descent set of $\sigma$.
Similarly we can define recoil compositions. 

If $I=(i_1,\dots,i_q)$ and $J=(j_1,\dots,j_p)$ are two compositions, then
$I \cdot J$ refers to their concatenation $(i_1,\dots,i_q,j_1,\dots,j_p)$,
and $I \mytr J$ is equal to $(i_1,\dots,i_q+j_1, j_2,\dots,j_p)$.

The {\it major index} $\maj(K)$ of a composition
$K = (k_1,\dots,k_r)$ is equal to the dot product of
$(k_1,\dots,k_r)$ with $(r-1,r-2,\dots,2,1,0)$, i.e.
$\sum_{i=1}^r (r-i) k_i$.

%%%%%%%%%%%%%%%%%%%%%%%%%%%%%%%%%%%%%%%%%%%%%%%%%%%%%%%%%%%%%%%%%%%%%%%%%%%%%%%
\subsection{Word Quasi-symmetric functions: $\WQSym$}

Let $w\in A^*$.
The \emph{packed word} $u=\pack(w)$ associated with $w$ is
obtained by the following process. If $b_1<b_2<\ldots <b_r$ are the letters
occuring in $w$, $u$ is the image of $w$ by the homomorphism
$b_i\mapsto a_i$.
A word $u$ is said to be \emph{packed} if $\pack(u)=u$. We denote by $\PW$ the
set of packed words.
With such a word, we associate the polynomial
\begin{equation}
\M_u :=\sum_{\pack(w)=u}w\,.
\end{equation}
{\footnotesize
For example, restricting $A$ to the first five integers,
\begin{equation}
%\begin{split}
\M_{13132}= 13132 + 14142 + 14143 + 24243 
+ 15152 + 15153 + 25253 + 15154 + 25254 + 35354.
%\end{split}
\end{equation}
}

Under the abelianization
$\chi:\ \K\langle A\rangle\rightarrow\K[X]$, the $\M_u$ are mapped to the
monomial quasi-symmetric functions $M_I$ 
($I=(|u|_a)_{a\in A}$ being the evaluation vector of $u$).

These polynomials span a subalgebra 
of $\K\langle A\rangle$, called $\WQSym$ for Word
Quasi-Symmetric functions~\cite{Hiv-adv}.
These are the invariants of the noncommutative
version of Hivert's quasi-symmetrizing action \cite{Hiv-adv}, which is
defined by
$\sigma\cdot w = w'$ where $w'$ is such that $\Std(w')=\Std(w)$ and
$\chi(w')=\sigma\cdot\chi(w)$.
Thus, two words are in the same $\SG(A)$-orbit iff they have the same packed
word.

The graded dimension of $\WQSym$ is the sequence of ordered Bell numbers
(\cite[A000670]{Slo}) $1, 1, 3, 13, 75, 541, 4683, 47293, 545835,\dots$.
Hence, $\WQSym$ is much larger than $\Sym$, which can be embedded in it in
various ways~\cite{NT4,NT07}.

The product of the $\M_u$ of $\WQSym$ is given by
\begin{equation}
\label{prodG-wq}
\M_{u'} \M_{u''} = \sum_{u \in u'\convW u''} \M_u\,,
\end{equation}
where the \emph{convolution} $u'\convW u''$ of two packed words
is defined as
\begin{equation}
u'\convW u'' = \sum_{v,w ;
u=v\cdot w\,\in\,\PW, \pack(v)=u', \pack(w)=u''} u\,.
\end{equation}

{\footnotesize
For example,
\begin{equation}
\label{M1121}
\M_{11} \M_{21} =
\M_{1121} + \M_{1132} + \M_{2221} + \M_{2231} + \M_{3321}.
\end{equation}
\begin{equation}
\begin{split}
\M_{21} \M_{121} =& \ \ \
\M_{12121} + \M_{12131} + \M_{12232} + \M_{12343} + \M_{13121} +
\M_{13232} + \M_{13242}\\
&+ \M_{14232} + \M_{23121} + \M_{23131} + \M_{23141} + \M_{24131} +
\M_{34121}.
\end{split}
\end{equation}
}

%%%%%%%%%%%%%%%%%%%%%%%%%%%%%%%%%%%%%%%%%%%%%%%%%%%%%%%
\subsection{Matrix quasi-symmetric functions: $\MQSym$}

This algebra is introduced in \cite{Hiv-adv, NCSF6}.
We start from a totally ordered set
of commutative variables $X=\{x_1<\cdots<x_n\}$ and consider the ideal
$\K[X]^+$ of polynomials without constant term.  We denote by
$\K\{X\}=T(\K[X]^+)$ its tensor algebra. 
We will also consider tensor products of elements of this
algebra. To avoid confusion, we denote by ``$\pP$'' the product
of the tensor algebra and call it the dot product. We reserve the
notation $\tensor$ for the external tensor product. 

A natural basis of $\K\{X\}$ is formed by dot products of
nonconstant monomials (called \emph{multiwords} in the sequel),
which can be represented by nonnegative integer matrices
$M=(m_{ij})$, where $m_{ij}$ is the exponent of the variable $x_i$
in the $j$th factor of the tensor product. Since constant monomials
are not allowed, such matrices have no zero column. We say that
they are \emph{horizontally packed}. 
A multiword $\mm$ can be encoded in the following
way. Let $V$ be the \emph{support} of $\mm$, that is, the set
of those variables $x_i$ such that the $i$th row of $M$ is non zero,
and let $P$ be the matrix obtained from $M$ by removing the null rows.
We set $\mm=V^P$. A matrix such as $P$, without zero rows or
columns, is said to be \emph{packed}.

{\footnotesize
For example the multiword $\mm = a\pP ab^3e^5\pP a^2d$ is
encoded by
$\indexmat\SMat{1&1&2\\0&3&0\\0&0&0\\0&0&1\\0&5&0}$.
Its support is the set $\{a,b,d,e\}$, and the associated packed matrix
is $\SMat{1&1&2\\0&3&0\\0&0&1\\0&5&0}$.
}

\bigskip
Let $\MQSym(X)$ be the linear subspace of $\K\{X\}$ spanned
by the elements
\begin{equation}
\MS_M = \sum_{V\in{\mathcal P}_k(X)}V^M
\end{equation}
where ${\mathcal P}_k(X)$ is the set of $k$-element subsets of $X$, and
$M$ runs over packed matrices of height $h(m)<n$.

\bigskip
{\footnotesize
For example, on the alphabet $\{a<b<c<d\}$
\renewcommand{\indexmat}%
    {\smallmatrice{\Hack a\\\Hack b\\\Hack c\\\Hack d\\}}
$$
     \MS_\SMat{1&1&2\\0&3&0\\0&0&1}=
     \indexmat\SMat{1&1&2\\0&3&0\\0&0&1\\0&0&0}+
     \indexmat\SMat{1&1&2\\0&3&0\\0&0&0\\0&0&1}+
     \indexmat\SMat{1&1&2\\0&0&0\\0&3&0\\0&0&1}+
     \indexmat\SMat{0&0&0\\1&1&2\\0&3&0\\0&0&1}
$$
}

One can show that $\MQSym$ is a subalgebra of $\K\{X\}$. Actually,
$$
\MS_P\MS_Q =\sum_{R\in\ash (P,Q)} \MS_R
$$
where the {\em augmented shuffle} of $P$ and $Q$,
$\ash(P,Q)$ is defined as follows: let $r$  be an
integer between $\max(p,q)$ and $p+q$, where $p=h(P)$ and $q=h(Q)$.
Insert null rows in the matrices $P$ and $Q$ so as to form 
matrices $\tilde P$ and $\tilde Q$ of height $r$. Let $R$ be
the matrix $(\tilde P,\tilde Q)$. The set $\ash (P,Q)$ is formed
by all the matrices without null rows obtained in this way.

{\footnotesize
For example :
$$
\begin{array}{l}
  \MS{\SMat{2&1\\1&0}}\MS_{\SMat{3&1}} = \\[3mm]
\qquad\ 
  \MS{\SMat{2&1&0&0\\1&0&0&0\\0&0&3&1}}+\MS{\SMat{2&1&0&0\\1&0&3&1}}+
  \MS{\SMat{2&1&0&0\\0&0&3&1\\1&0&0&0}}+\MS{\SMat{2&1&3&1\\1&0&0&0}}+
  \MS{\SMat{0&0&3&1\\2&1&0&0\\1&0&0&0}}
\end{array}
$$
}

%%%%%%%%%%%%%%%%%%%%%%%%%%%%%%%%%%%%%%%%%%%%%%%%%%%%%%%%%%%%%%%%%%%%%%%%%%%%%%%
\subsection{Free quasi-symmetric functions: $\FQSym$}

The Hopf algebra $\FQSym$ is the subalgebra of $\WQSym$ spanned by the
polynomials~\cite{NCSF6}
\begin{equation}
\G_\sigma := \sum_{\Std(u)=\sigma} \M_u
= \sum_{\Std(w)=\sigma} w.
\end{equation}
The multiplication rule is, for $\alpha\in\SG_k$ and
$\beta\in\SG_l$,
\begin{equation}\label{multG}
\G_\alpha \G_\beta
= \sum_{\genfrac{}{}{0pt}{}{\gamma\in\SG_{k+l};\,\gamma=u\cdot v}
  {\Std(u)=\alpha,\Std(v)=\beta}}\G_\gamma\,.
\end{equation}
As a Hopf algebra, $\FQSym$ is self-dual. The scalar
product materializing this duality is the one for which
$(\G_\sigma\,,\,\G_\tau)=\delta_{\sigma,\tau^{-1}}$
(Kronecker symbol).
Hence, $\F_\sigma:=\G_{\sigma^{-1}}$ is the dual basis of $\G$.
Their product is given by 
\begin{equation}
\label{multF}
\F_\alpha \F_\beta = \sum_{\gamma\in\alpha\ssh\beta} \F_\gamma.
\end{equation}

%%%%%%%%%%%%%%%%%%%%%%%%%%%%%%%%%%%%%%%%%%%%%%%%%%%%%%%%%%%%%%%%%%%%%%%%%%%%%%%
\subsection{Embeddings}

%%%%%%%%%%%%%%%%%%%%%%%%%%%%%%%%%%%%
\subsubsection{$\Sym$ into $\MQSym$}

Recall that the algebra of {\it noncommutative symmetric functions}
is the free associative algebra
$\Sym=\C \langle S_1,S_2,\dots \rangle$ generated by an infinite sequence of
non-commutative indeterminates $S_k$, called \emph{complete} symmetric
functions.
For a composition $I=(i_1,\dots,i_r)$, one sets $S^I = S_{i_1} \dots S_{i_r}$.
The family $(S^I)$ is a linear basis of $\Sym$.
A useful realization, denoted by $\Sym(A)$, can be obtained by taking an
infinite alphabet $A=\{a_1,a_2,\dots \}$ and defining its complete homogeneous
symmetric functions by the generating function
\begin{equation*}
\sum_{n \geq 0} t^n S_n(A) = (1-ta_1)^{-1}(1-ta_2)^{-1} \dots .
\end{equation*}

Given a packed matrix $P$, the vector of its \emph{column sums} will be
denoted by $\text{Col}(P)$.
The algebra morphism defined on generators by
\begin {equation}
\beta:\ S_n\mapsto \sum_{{\rm Col}(P)=(n)}\MS_P
\end{equation}
is an embedding of Hopf algebras \cite{NCSF6}.
By definition of $\MQSym$, for an arbitrary composition, we have
\begin{equation}
\beta(S^I) = \sum_{{\rm Col}(P)=I}\MS_P .
\end{equation}

\subsubsection{$\Sym$ into $\WQSym$}

The algebra morphism defined on generators by
\begin {equation}
\alpha:\ S_n\mapsto \sum_{\Std(u)=12\cdots n}\M_u
\end{equation}
is also an embedding of Hopf algebras. For an arbitrary composition,
\begin{equation}
\alpha(S^I) = \sum_{\DesC(u)\fatter I}\M_u\,.
\end{equation}
Indeed, when $\Sym$ is realized as $\Sym(A)$, the latter
sum is equal to $S^I(A)$.

%%%%%%%%%%%%%%%%%%%%%%%%%%%%%%%%%%%%%%%%%%%%%%%%%%%%%%%%%%%%%%%%%%%%%%%%%%%%%%%
\subsection{Epimorphisms}

We shall also need to project back from the algebras $\MQSym$ and $\WQSym$ to
$\Sym$.
The crucial projection is the one associated with the 
(non-Hopf) quotient of $\WQSym$ introduced
in \cite{HNTT}.

\subsubsection{$\MQSym$ to $\WQSym$}

To a packed matrix $M$, one associates a packed word $w(M)$
as follows. Read the entries of $M$ columnwise, from top to
bottom and left to right. The word $w(M)$ is obtained by
repeating $m_{ij}$ times each row index $i$.

Let ${\mathcal J}$ be the ideal of $\MQSym$ generated
by the differences 
\begin{equation}
\{\MS_P-\MS_Q |w(P)=w(Q)\}.
\end{equation}
Then the quotient $\MQSym/{\mathcal J}$ is isomorphic
as an algebra to $\WQSym$,
via the identification $\overline{\MS}_M=\M_{w(M)}$.
More precisely, 
$\eta:\ \overline{\MS}_M\mapsto \M_{w(M)}$ is a morphism of algebras.

\subsubsection{$\WQSym$ to $\Sym$}

Let $w$ be a packed word.
The \emph{Word composition} (W-comp\-os\-ition) of $w$ is the composition
whose descent set is given by the positions of the last occurrences of each
letter in $w$.

{\footnotesize
For example,
\begin{equation}
\WC(1543421323) = (2,3,2,2,1).
\end{equation}
Indeed, the descent set is $\{2,5,7,9,10\}$ since
the last $5$ is in position $2$,
the last $4$ is in position $5$,
the last $1$ is in position $7$,
the last $2$ is in position $9$,
and the last $3$ is in position $10$.

The following tables group the packed words in
$\PW_2$ and $\PW_3$ according to their W-composition.

\begin{equation}
\begin{array}{|c|c|}
\hline
2 & 11 \\
\hline
\hline
11 & 12 \\
   & 21 \\
\hline
\end{array}
\hskip2cm
\begin{array}{|c|c|c|c|}
\hline
3   & 21  & 12  & 111\\
\hline
\hline
111 & 112 & 122 & 123\\
    & 121 & 211 & 132\\
    & 212 &     & 213\\
    & 221 &     & 231\\
    &     &     & 312\\
    &     &     & 321\\
\hline
\end{array}
\end{equation}
}

%%%%%%%%%%%%%%%%%%%%%%%%%%%%%%%%%%%%%%%%%%%%%%%%%%%%%%%%%%%%%%%%%%%%%%%%%%%%%%%
Let $\sim$ be the equivalence relation on packed words defined by
$u\sim v$ iff $\WC(u)=\WC(v)$.
Let ${\mathcal J'}$ be the subspace of $\WQSym$ spanned by the differences
\begin{equation}
\{ \M_u - \M_v\, |\, u\sim v\}.
\end{equation}

Then, it has been shown \cite{HNTT} that ${\mathcal J'}$ is a two-sided
ideal of $\WQSym$, and that the quotient $\TEb$ defined by
$\TEb=\WQSym/{\mathcal J'}$ is isomorphic to $\NCSF$ as an algebra.

More precisely, recall that $\Psi_n$ is a {\it noncommutative power sum of the
first kind}~\cite{NCSF1}.
Tevlin defined the {\it noncommutative monomial symmetric functions} $\Psi_I$
\cite{Tev} as quasideterminants in the $\Psi_n$'s.
We do not need the precise definition of $\Psi_I$ here, only the following
result.

\begin{proposition}\cite{HNTT}\label{HNTTMorphism}
$\zeta:\ \overline{\M}_u\mapsto \Psi_{\WC(u)}$ is a morphism
of algebras.
\end{proposition}

%%%%%%%%%%%%%%%%%%%%%%%%%%%%%%%%%%%%%%%%%%%%%%%%%%%%%%%%%%%%%%%%%%%%%%%%%%%%%%%
%%%%%%%%%%%%%%%%%%%%%%%%%%%%%%%%%%%%%%%%%%%%%%%%%%%%%%%%%%%%%%%%%%%%%%%%%%%%%%%
%%%%%%%%%%%%%%%%%%%%%%%%%%%%%%%%%%%%%%%%%%%%%%%%%%%%%%%%%%%%%%%%%%%%%%%%%%%%%%%
\section{Quantizations and  noncommutative Hall-Littlewood functions} 

In this section, we introduce a new $q$-analogue $S^I(q)$ of the basis $S^I$
of $\Sym$, giving two different but equivalent definitions.
When we examine the transition matrices between this new basis and other
bases, we will see a connection to permutation tableaux and hence
to the asymmetric exclusion process.
The new basis elements $S^I(q)$ play the role of the classical Hall-Littlewood
$Q'_\mu$~\cite[Ex. 7.(a) p. 234]{Mcd}, and of Hivert's $H_I(q)$.

%%%%%%%%%%%%%%%%%%%%%%%%%%%%%%%%%%%%%%%%%%%%%%%%%%%%%%%%%%%%%%%%%%%%%%%%%%%%%%%
\subsection{The special inversion statistic}

Let $u=u_1\cdots u_n$ be a packed word. We say that an inversion $u_i=b>u_j=a$
(where $i<j$ and $a<b$) is {\em special} if $u_j$ is the {\em rightmost}
occurence of $a$ in $u$. Let $\sinv (u)$ denote the number of special
inversions in $u$.
Note that if $u$ is a permutation, this coincides with its ordinary inversion
number.

%%%%%%%%%%%%%%%%%%%%%%%%%%%%%%%%%%%%%%%%%%%%%%%%%%%%%%%%%%%%%%%%%%%%%%%%%%%%%%%
\subsection{Quantizing $\WQSym$}

Let $\M'_u=q^{\sinv(u)}\M_u$ and define a linear map $\phi_q$ by
$\phi_q(\M_u)=\M'_u$.
We define a new associative product $\qp$ on $\WQSym$ by requiring that
\begin{equation}
\M'_u\qp \M'_v=\phi_q(\M_u\M_v)\,.
\end{equation}

{\footnotesize
For example, by~(\ref{M1121}), one has
\begin{equation}
\begin{split}
\M'_{11} \qp \M'_{21}
&= \M'_{1121} + \M'_{1132} + \M'_{2221} + \M'_{2231} + \M'_{3321}\\
&= q\M_{1121} + q\M_{1132} + q^3\M_{2221} + q^3\M_{2231} + q^5\M_{3321}.
\end{split}
\end{equation}
}

This algebra structure on the vector space $\WQSym$
will be denoted by $\WQSym_q$.

%%%%%%%%%%%%%%%%%%%%%%%%%%%%%%%%%%%%%%%%%%%%%%%%%%%%%%%%%%%%%%%%%%%%%%%%%%%%%%%
\subsection{Quantizing $\MQSym$}

Similarly, the $q$-product $\qp$ can be defined on $\MQSym$,
by requiring that the $\MS'_M=q^{\sinv(w(M))}\MS_M$ multiply as
the $\MS_M$.

%%%%%%%%%%%%%%%%%%%%%%%%%%%%%%%%%%%%%%%%%%%%%%%%%%%%%%%%%%%%%%%%%%%%%%%%%%%%%%%
%%%%%%%%%%%%%%%%%%%%%%%%%%%%%%%%%%%%%%%%%%%%%%%%%%%%%%%%%%%%%%%%%%%%%%%%%%%%%%%
%%%%%%%%%%%%%%%%%%%%%%%%%%%%%%%%%%%%%%%%%%%%%%%%%%%%%%%%%%%%%%%%%%%%%%%%%%%%%%%
%%%%%%%%%%%%%%%%%%%%%%%%%%%%%%%%%%%%%%%%%%%%%%%%%%%%%%%%%%%%%%%%%%%%%%%%%%%%%%%
\subsection{Two equivalent definitions of $S^I(q)$}

Embedding $\Sym$ into $\MQSym_q$ and projecting back to $\Sym$,
we define $q$-analogues of the products $S^I$ by
\begin{equation}
S^I(q)=\zeta\circ\eta ( \beta(S_{i_1})\qp\cdots\qp \beta(S_{i_r}))\,.
\end{equation}

\noindent
Equivalently, since under the above embeddings, the image of $\Sym$ in
$\MQSym$ is contained in the image of $\WQSym$, one can embed $\Sym$ into
$\WQSym_q$ and project back to $\Sym$, which yields
\begin{equation}
S^I(q)=\zeta (\alpha(S_{i_1})\qp\cdots\qp \alpha(S_{i_r}))\,.
\end{equation}

%%%%%%%%%%%%%%%%%%%%%%%%%%%%%%%%%%%%%%%%%%%%%%%%%%%%%%%%%%%%%%%%%%%%%%%%%%%%%%%
\subsection{The transition matrix $M(S(q),\Psi)$}

For any two bases $F$, $G$ of $\Sym$, we denote by $M_n(F,G)$
the matrix indexed by compositions of $n$, whose entry in
row $I$ and column $J$ is the coefficient of $G_I$ in the
$G$-expansion of $F_J$. 
We will give two 
combinatorial formulas (Propositions \ref{SP1} and \ref{SP2})
and one recursive formula (Theorem \ref{recursive})
for the elements of the transition matrix 
$M(S(q),\Psi)$, where the $\Psi_I$'s are Tevlin's noncommutative
monomial symmetric functions.

%%%%%%%%%%%%%%%%%%%%%%%%%%%%%%
\subsubsection{First examples}

Let $[n]$ denote the $q$-analogue $1+q+\dots+q^{n-1}$ of $n$.
The first transition matrices $SP_n=M_n(S(q),\Psi)$ are 

\begin{equation*}
SP_3 = M_3(S(q),\Psi) =
\left(
\begin{matrix}
[1] & [1] & [1] & [1]    \\
[1] & [3] & [2] & [2][2] \\
[1] & [1] & [2] & [2]    \\
[1] & [3] & [3] & [2][3] 
\end{matrix}
\right)
\end{equation*}

\begin{equation*}
SP_4 =
\left(
\begin{matrix}
[1] & [1] & [1]        & [1]    & [1] & [1]    & [1]    & [1]       \\
[1] & [4] & [3]        & [2][3] & [2] & [2][3] & [2][2] & [2][2][2] \\
[1] & [1] & [3]        & [3]    & [2] & [2]    & [2][2] & [2][2]    \\
[1] & [4] & [4][3]/[2] & [3][4] & [3] & [3][3] & [3][3] & [2][3][3] \\
[1] & [1] & [1]        & [1]    & [2] & [2]    & [2]    & [2]       \\
[1] & [4] & [3]        & [2][3] & [3] & [3][3] & [2][3] & [2][2][3] \\
[1] & [1] & [3]        & [3]    & [3] & [3]    & [2][3] & [2][3]    \\
[1] & [4] & [4][3]/[2] & [3][4] & [4] & [3][4] & [3][4] & [2][3][4]
\end{matrix}
\right)
\end{equation*}

\bigskip
The coefficient of $\Psi_I$ in $S^J(q)$ will be denoted by $C_I^J(q)$.

%%%%%%%%%%%%%%%%%%%%%%%%%%%%%%%%%%%%%%%%%%%%%
\subsubsection{Combinatorial interpretations}

Recall that by Proposition \ref{HNTTMorphism}, 
$\zeta\circ\eta$ is a morphism of algebras
sending $\MS_M$ to $\Psi_{\WC(w(M))}$. Hence, our first 
definition of $S^J(q)$ gives the following:

\begin{proposition}\label{SP1}
Let $I=(i_1,\dots,i_k)$ and $J=(j_1,\dots,j_l)$ be two compositions.
Let $M(I,J)$ be the set of integer matrices
$M=(m_{p,q})_{1\leq p\leq l;1\leq q\leq k}$ without null rows such that
\begin{equation}
\WC(w(M))=I \qquad\text{and}\qquad \Col(M)=J.
\end{equation}
Then
\begin{equation}
C_I^J(q) = \sum_{M\in M(I,J)} q^{\sinv(w(M))}.
\end{equation}
\end{proposition}

{\footnotesize
For example, the six matrices corresponding to the coefficient $[4][3]/[2]$ of
$M_4$ in row $(2,1,1)$ and column $(2,2)$ are
\begin{equation}
\label{mats211}
\SMat{2&.\\.&1\\.&1}
\quad
\SMat{1&1\\1&.\\.&1}
\quad
\SMat{1&1\\.&1\\1&.}
\quad
\SMat{.&1\\2&.\\.&1}
\quad
\SMat{.&1\\1&1\\1&.}
\quad
\SMat{.&1\\.&1\\2&.}
\end{equation}
The corresponding statistics are
\begin{equation}
\{ 0,\ 1,\ 2,\ 2,\ 3,\ 4 \}.
\end{equation}
}

The second definition of $S^J(q)$ yields a different combinatorial description:

\begin{proposition}
\label{SP2}
\label{cijW}
Let $I$ and $J$ be two compositions and let $W(I,J)$ be the set of packed
words $w$ such that
\begin{equation}
\WC(w)=I \qquad\text{and}\qquad \DesC(w)\fatter J.
\end{equation}
Then
\begin{equation}
C_I^J(q) = \sum_{w\in W(I,J)} q^{\sinv(w)}.
\end{equation}
\end{proposition}

{\footnotesize
For example, the six packed words corresponding to the coefficient
$[4][3]/[2]$ of $M_4$ in row $I=(2,1,1)$ and column $J=(2,2)$ are
\begin{equation}
1123,\ 1213,\ 1312,\ 2213,\ 2312,\ 3312.
\end{equation}
These words are the column readings of the six matrices 
from (\ref{mats211}).  The first word has descent composition
$(4)$ and the others have descent composition $(2,2)$.
}

%%%%%%%%%%%%%%%%%%%%%%%%%%%%%%%%%%%%%%%%%%
\subsubsection{The $q$-product on $\NCSF$}

To explain the factorization of  the coefficients of the transition
matrix $M(S(q),\Psi)$, we  need a recursive formula for $S^J(q)$.
This will be given in Theorem \ref{recursive}.

In $\WQSym$, let 
\begin{equation}
\tS_n=\alpha(S_n)=\sum_{u;u\uparrow n}\M_u
\end{equation}
where $u\uparrow n$ means that $u$ is a nondecreasing packed word of length
$n$, and define
\begin{equation}
\tS^J=\tS_{j_1}\qp\dots\qp \tS_{j_r}
\end{equation}
so that $S^J(q)=\zeta(\tS^J)$.
Let  $J=(j_1,\dots,j_r)$ and set $J'=(j_2,\dots,j_r)$.
Since $\qp$ is associative in $\WQSym_q$, we have 
\begin{equation}
S^{j_1,J'}(q) =
\zeta( \tS^{j_1,J'}) =
\sum_{\genfrac{}{}{0pt}{}{u,v;u\uparrow j_1}{\DesC(v)\fatter J'}}
      q^{\sinv(v)} \zeta( \M_u \qp \M_v).
\end{equation}
This expression can be simplified by means of the following Lemma.
\begin{lemma}
\label{lem-uvuvp}
Let $u$ be a nondecreasing packed word. Then
\begin{equation}
\zeta(\M_u\qp\M_v)=\zeta(\M_u\qp\M_{v'})
\end{equation}
for all $v'$ such that $\WC(v')=\WC(v)$.
\end{lemma}

\Proof
Since $u$ is nondecreasing, each packed word  $z=x\cdot y$ appearing in the
expansion of $\M_u\qp\M_v$ is completely determined by the letters used in $x$
and the letters used in $y$.
Looking at the packed words $z$ and $z'$ occuring in  $\M_u\qp\M_v$ and in
$\M_u\qp\M_{v'}$ with given letters used for their prefixes and suffixes, we
have $\sinv(z')=\sinv(z)+\sinv(v')-\sinv(v)$, whence the result.
\qed

{\footnotesize
For example,
\begin{equation}
\M_{11}\qp \M_{12} = \M_{1112} + \M_{1123} + q^2\M_{2212}
                   + q^2\M_{2213} + q^4\M_{3312}.
\end{equation}
\begin{equation}
\M_{11}\qp q\M_{21} = q\M_{1121} + q\M_{1132} + q^3\M_{2221}
                   + q^3\M_{2231} + q^5\M_{3321}.
\end{equation}
}

Let now $\sigma:\NCSF\to\WQSym$ be the section of the projection $\zeta$
defined by
\begin{equation}
\sigma(\Psi_I) = \M_{1^{i_1}2^{i_2}\dots r^{i_r}}.
\end{equation}
We can define a (non-associative!) $q$-product on $\NCSF$ by
\begin{equation}\label{non-assoc}
f \qp g = \zeta( \sigma(f) \qp \sigma(g)).
\end{equation}
Then Lemma~\ref{lem-uvuvp} implies that
\begin{equation}
S^I(q) = S^{i_1} \qp (S^{i_2}\qp(\dots (S^{i_{r-1}}\qp S^{i_r}))).
\end{equation}

%%%%%%%%%%%%%%%%%%%%%%%%%%%%%%%%%%%%%%%%%%%%%%%%
\subsubsection{Closed form for the coefficients}

>From Lemma~\ref{lem-uvuvp}, we now have
\begin{equation}
S^{j_1,J'}(q) =
\zeta( \tS^{j_1,J'}) =
\sum_{\genfrac{}{}{0pt}{}{u,v;u\uparrow j_1}{v\uparrow j_2+\dots+j_r}}
      C^{J'}_{\WC(v)}(q) \zeta( \M_u \qp \M_v).
\end{equation}
Note that $\zeta(\M_u\qp\M_v)$ and $\zeta(\M_{u'}\qp\M_{v'})$ are linear
combinations of
disjoint sets of $\Psi_K$ as soon as the nondecreasing words $v$ and $v'$ are
different. So the computation of the coefficient $C_J^I$ boils down to
the evaluation of
\begin{equation}
\sum_{u;u\uparrow j_1} \zeta( \M_u \qp \M_v) = \zeta(\tS^{j_1}\qp\M_v),
\end{equation}
where $v$ is a nondecreasing word.
Let us first  characterize the terms of the product yielding a given $\Psi_I$.
\begin{lemma}
\label{QR}
Let $u$ be a nondecreasing word of length $k$ over $[1,r]$.
Given a composition $I=(i_1,\dots,i_r)$ of length $r$,
there exists at most one nondecreasing word $v$ over $[1,r]$ such that
$uv$ is packed and $\WC(uv)=I$.
Such a $v$ exists precisely when $u=u_1\cdots u_k$ satisfies $u_i<u_{i+1}$ for
$i\in\Des(I)$.

In this case, let $y=1^{i_1}2^{i_2}\dots r^{i_r}$.
Then $\sinv(uv)$ is equal to
\begin{equation}
\sum_{1\leq i\leq k} (u_i-y_i).
\end{equation}
This sum is also equal to 
\begin{equation}
\sum_{1\leq i\leq k} u_i -  (k+\maj(\overline K)),
\end{equation}
where $K$ is the composition of $k$ such that $\Des(K)=\Des(I)\cap[1,k-1]$.
\end{lemma}

\Proof
The construction of $v$ was already given in the proof of Theorem 6.1
of~\cite{HNTT}. It comes essentially from the facts that the letters which
should be used in $v$ are determined by the letters used in $u$, and that a word is
uniquely determined  by its packed word and its alphabet.

Now, for each letter $x$ of $u$, its contribution to $\sinv(uv)$ is given by
the number of different letters strictly smaller than $x$ appearing in $v$.
This is equal to $u_i-y_i$. The sum of the $y_i$ is $k+\maj(\overline K)$.
\qed

{\footnotesize
For example, given $I=1221$, there are 10 nondecreasing words $u$ of $[1,4]$
of length $3$ satisfying the conditions of the lemma. The following table
gives the corresponding~$v$ and the $\sinv$ statistics of the products $uv$.
\begin{equation}
\label{tab1221}
\begin{array}{|c|c|c|}
\hline
u   & v   & \sinv \\
\hline
\hline
122 & 334 & 0   \\
123 & 224 & 1   \\
124 & 223 & 2   \\
133 & 224 & 2   \\
134 & 223 & 3   \\
144 & 223 & 4   \\
233 & 114 & 3   \\
234 & 113 & 4   \\
244 & 113 & 5   \\
344 & 112 & 6   \\
\hline
\end{array}
\end{equation}
}

We are now in a position to compute
\begin{equation}
\zeta(\tS^{j_1}\qp\M_v)
\end{equation}
when $v$ is a nondecreasing word.
\begin{lemma}
\label{qprod-init}
Let $I$ be a composition of $k+n$ and let $I'$ be the composition of $n$ such
that $\Des(I')=\{a_1,\dots,a_s\}$ satisfies
$\{k+a_1,\dots,k+a_s\} = \Des(I)\cap[k+1,k+n]$.

Let $v$ be the nondecreasing word of evaluation $I'$.
The coefficient of $\Psi_I$ in $\zeta(\tS_k\qp\M_v)$
is the $q$-binomial coefficient
\begin{equation}
\qbin{k+r-s}{r-s}
\end{equation}
where $r=l(I)$, $K$ is the composition of $k$ such that
$\Des(K)=\Des(I)\cap[1,k-1]$, and $s=l(K)$.
\end{lemma}

\Proof
To start with, write
\begin{equation}
\tS_k\qp\M_v=\sum_w q^{\sinv(w)}\M_w,
\end{equation}
where $w$ runs over packed words of the form $w=u'v'$, with
$u'$ nondecreasing and $\pack(v')=v$. From Lemma \ref{QR}, we
see that in order to have $\WC(u'v')=I$, $u'=x_1\cdots x_k$
must be a word over the interval $[1,r]$ with equalities $x_i=x_j$
allowed precisely when cells $i$ and $j$ are in the same row
of the diagram of $K$. The commutative image of the formal
sum of such words, which are the nondecreasing reorderings of
the quasi-ribbons of shape $K$ \cite{NCSF4}, is the quasi-symmetric
{\it quasi-ribbon } polynomial $F_K(t_1,\ldots,t_r)$,
introduced in \cite{Gessel}.  Hence, the coefficient of
$\Psi_I$ is
\begin{equation}
q^{-\maj(\bar K)}\qbin{k+r-s}{r-s}
\end{equation} 
given the generating function~\cite{GR}
\begin{equation}
\sum_{m\ge 0}t^mF_K(1,q,\ldots,q^{m-1})
     =\frac{t^{l(K)} q^{\maj(\bar K)}}{(t;q)_{k+1}}\,.
\end{equation}
\qed

{\footnotesize
The example presented in~(\ref{tab1221}) corresponds to the case $I=(1,2,2,1)$
and $i=3$, so that $K=(1,2)$. We then find $\qbin{3+4-2}{4-2}=\qbin{5}{2}$,
which indeed corresponds to the statistic in the last column
of~(\ref{tab1221}).
}

Summarizing the above discussion, we can now state the main result of
this section:

\begin{theorem}
\label{recursive}
Let $I=(i_1,\dots,i_k)$ and $J=(j_1,\dots,j_l)$ be compositions of $n$.
Then the coefficient $C_I^J(q)$ of $\Psi_I$ in $S^J(q)$ is given by the
following rule:\\
(i) if $i_1<j_1$, then $C_I^J(q)=C_{(i_1+i_2,i_3,\dots,i_k)}^J(q)$,\\
(ii) otherwise,
\begin{equation}
C_I^J(q) = \qbin{k+j_1-1}{j_1} C_{I'}^{(j_2,\dots,j_l)}(q)
\end{equation}
where the diagram of $I'$ is obtained by removing the
first $j_1$ cells of the diagram of $I$.
\end{theorem}
\hfill
\qed

%%%%%%%%%%%%%%%%%%%%%%%%%%%%%%%%%%%%%%%%%%%%%%%%%%%%%%%%%%%%%%%%%%%%%%%%%%%%%%%
\subsection{The transition matrix $M(R(q),\Psi)$}

%%%%%%%%%%%%%%%%%%%%%%%%%%%%%%%%%%%%
\subsubsection{$q$-deformed ribbons}

We now define a $q$-ribbon basis $R_I(q)$ in terms of the 
$S^J(q)$'s by analogy to the relationship between 
the ordinary $R_I$'s and $S^J$'s:
\begin{equation}
R_I(q) := \sum_{J\fatter I} (-1)^{l(J)-l(I)} S^J(q).
\end{equation}
The coefficient of $\Psi_I$ in the expansion of $R^J(q)$ will be denoted
by $D_I^J(q)$.

%%%%%%%%%%%%%%%%%%%%%%%%%%%%%%
\subsubsection{First examples}

We get the following transition matrices between $R(q)$ and $\Psi$ 
for $n=3,4$:

\begin{equation*}
RP_3 = M_3(R,\Psi) =
\left(
\begin{matrix}
1 & . & . & . \\
1 & q+q^2 & q & . \\
1 & . & q & . \\
1 & q+q^2 & q+q^2 & q^3 
\end{matrix}
\right)
\end{equation*}

\begin{equation*}
RP_4 = 
\left(
\begin{matrix}
1 & . & . & . & . & . & . & . \\
1 & q[3] & q\!+\!q^2 & . & q & q^2 & . & . \\
1 & . & q\!+\!q^2 & . & q & . & . & . \\
1 & q[3] & q\!+\!2q^2\!+\!q^3\!+\!q^4 & q^3[3]%\!+\!q^4\!+\!q^5
  & q\!+\!q^2 & q^2\!+\!q^3\!+\!q^4 & q^3 & . \\
1 & . & . & . & q & . & . & . \\
1 & q[3] & q\!+\!q^2 & . & q\!+\!q^2 & q^2\!+\!q^3\!+\!q^4
  & q^3 & . \\
1 & . & q\!+\!q^2 & . & q\!+\!q^2 & . & q^3 & . \\
1 & q[3] & q\!+\!2q^2\!+\!q^3\!+\!q^4 & q^3[3]%q^3\!+\!q^4\!+\!q^5
  & q[3] & q^2\!+\!q^3\!+\!2q^4\!+\!q^5 & q^3[3] & q^6
\end{matrix}
\right)
\end{equation*}

%%%%%%%%%%%%%%%%%%%%%%%%%%%%%%%%%%%%%%%%%%%%%
\subsubsection{Combinatorial interpretations}

By definition of the transition matrix from $S(q)$ to $R(q)$,
the matrices $M(R(q),\Psi)$ can be described as follows:

\begin{proposition}
\label{dijW}
Let $I$ and $J$ be compositions of $n$, and let $W'(I,J)$ be the set of packed words
$w$ such that
\begin{equation}
\WC(w)=I \qquad\text{and}\qquad \DesC(w)=J.
\end{equation}
Then
\begin{equation}
D_I^J(q) = \sum_{w\in W'(I,J)} q^{\sinv(w)}.
\end{equation}
\end{proposition}

\Proof
This follows directly from the combinatorial interpretation of $C_I^J$
in terms of packed words (see Proposition~\ref{cijW}).
\qed

In terms of $\MQSym$, this can be rewritten as follows:
\begin{corollary}
Let $I$ and $J$ be compositions of $n$. Then $D_I^J(q)$ 
is given by the statistic
$\sinv(w(M))$ applied to the elements $M$ of the subset of $M(I,J)$ where in
each pair of consecutive columns, the bottommost nonzero entry of the left one
is strictly below the top-most nonzero entry of the right one.
\end{corollary}

%%%%%%%%%%%%%%%%%%%%%%%%%%%%%%%%%%%%%%%%%%%%%%%%%%%%%%%%%%%%%%%%%%%%%%%%%%%%%%%
\subsection{The transition matrix $M(L(q),\Psi)$}

%%%%%%%%%%%%%%%%%%%%%%%%%%%%%%%%%%%%%%%%%%%%%%%%%%%%%%%%%%%%
\subsubsection{A new $q$-analogue of the $L$ basis of $\NCSF$}

Let $\st(I,J)$ be the statistic on pairs of compositions of the same weight
defined by
\begin{equation}
\st(I,J):=
\left\{
\begin{array}{ll}
\#\{(i,j)\in \Des(I)\times\Des(J) | i\geq j\} & \text{if } I\finer J,\\
-\infty & \text{otherwise}
\end{array}
\right.
\end{equation}

We define a new basis $L(q)$  by
\begin{equation}
L_J(q) := \sum_{I\models |J|} q^{st(I,J)} \Psi_I
        = \sum_{I\finer J} q^{st(I,J)} \Psi_I.
\end{equation}
For $q=1$, this reduces to Tevlin's basis $L_I$ (in the notation
of \cite{HNTT}).
Since $M(L_J(q),\Psi)$ is unitriangular,
$L_I(q)$ is a basis of $\NCSF$.

%%%%%%%%%%%%%%%%%%%%%%%%%%%%%%
\subsubsection{First examples}

Here are the first transition matrices from $L(q)$ to $\Psi$:

\begin{equation} 
MLP_3 = M_3(L(q),\Psi) =
\left(
\begin{matrix}
1 & . & .   & . \\
1 & q & .   & . \\
1 & . & q   & . \\
1 & q & q^2 & q^3 
\end{matrix}
\right)
\end{equation}

\begin{equation}
MLP_4 =
%M_4(L(q),\Psi) = 
\left(
\begin{matrix}
1 & . & .   & .   & .   & .   & .   & .   \\
1 & q & .   & .   & .   & .   & .   & .   \\
1 & . & q   & .   & .   & .   & .   & .   \\
1 & q & q^2 & q^3 & .   & .   & .   & .   \\
1 & . & .   & .   & q   & .   & .   & .   \\
1 & q & .   & .   & q^2 & q^3 & .   & .   \\
1 & . & q   & .   & q^2 & .   & q^3 & .   \\
1 & q & q^2 & q^3 & q^3 & q^4 & q^5 & q^6
\end{matrix}
\right)
\end{equation}

Note that up to some minor changes (conjugation w.r.t.  mirror image of
compositions), these are the matrices expressing Hivert's Hall-Littlewood
$\tilde{H}_J$ on the basis $R_I$~\cite{Hiv-adv}.
This allows us to derive the expression of their inverse, that is, transition
matrices from $\Psi$ to $L(q)$ (see~\cite{Hiv-adv}, Theorem 6.6):
\begin{equation}
\label{PsiL}
\Psi_J = \sum_{I\finer J} (-1/q)^{l(I)-l(J)} q^{-\stb(I,J)} L_I(q),
\end{equation}
where $\stb(I,J)$ is
\begin{equation}
\stb(I,J):=
\left\{
\begin{array}{ll}
\#\{(i,j)\in \Des(I)\times\Des(J) | i\leq j\} & \text{if } I\finer J,\\
-\infty & \text{otherwise}
\end{array}
\right.
\end{equation}

%%%%%%%%%%%%%%%%%%%%%%%%%%%%%%%%%%%%%%%%%%%%%%%%%%%%%%%%%%%%%%%%%%%%%%%%%%%%%%%
\subsection{The transition matrix $M(S(q),L(q))$}
\label{sec-SL}

The coefficient of $L_I(q)$ in  $S^J(q)$ will be denoted by $E_I^J(q)$. In
this section, we will see a connection to permutation tableaux
and hence to the asymmetric exclusion process.

%%%%%%%%%%%%%%%%%%%%%%%%%%%%%%
\subsubsection{First examples}
\label{FirstExamples}

Here are the first transition matrices from $S(q)$ to $L(q)$:

\begin{equation*}
SL_3 = M_3(S(q),L(q)) =
\left(
\begin{matrix}
1 & 1 & 1 & 1 \\
. & 1+q & 1 & 2+q \\
. & . & 1 & 1 \\
. & . & . & 1 
%1 & 1 & 1 & 1 \\
%. & q+q^2 & q & 2q+q^2 \\
%. & . & q & q \\
%. & . & . & q^3 
\end{matrix}
\right)
\end{equation*}

\begin{equation*}
SL_4=
\left(
\begin{matrix}
1 & 1 & 1 & 1 & 1 & 1 & 1 & 1 \\
. & 1\!+\!q\!+\!q^2 & 1+q & 2+2q+q^2 & 1 & 2\!+\!2q\!+\!q^2 & 2+q & 3+3q+q^2 \\
. & . & 1+q & 1+q & 1 & 1 & 2+q & 2+q \\
. & . & q & 1+2q+q^2 & . & 1+q & 1+q & 3+3q+q^2 \\
. & . & . & . & 1 & 1 & 1 & 1 \\
. & . & . & . & . & 1+q & 1 & 2+q \\
. & . & . & . & . & . & 1 & 1 \\
. & . & . & . & . & . & . & 1
%1 & 1 & 1 & 1 & 1 & 1 & 1 & 1 \\
%. & q\!+\!q^2\!+\!q^3 & q\!+\!q^2 & 2q\!+\!2q^2\!+\!q^3 & q
%  & 2q\!+\!2q^2\!+\!q^3 & 2q\!+\!q^2 & 3q\!+\!3q^2\!+\!q^3 \\
%. & . & q\!+\!q^2 & q\!+\!q^2 & q & q & 2q\!+\!q^2 & 2q\!+\!q^2 \\
%. & . & q^4 & q^3\!+\!2q^4\!+\!q^5 & . & q^3\!+\!q^4 & q^3\!+\!q^4
%  & 3q^3\!+\!3q^4\!+\!q^5 \\
%. & . & . & . & q & q & q & q \\
%. & . & . & . & . & q^3+q^4 & q^3 & 2q^3\!+\!q^4 \\
%. & . & . & . & . & . & q^3 & q^3 \\
%. & . & . & . & . & . & . & q^6
\end{matrix}
\right)
\end{equation*}

In fact the right-hand column of each of these matrices contains
the (un-normalized) steady-state probabilities of each state of the 
partially asymmetric exclusion process (PASEP).  More specifically, 
the steady-state probabilities of the states 
$\bullet \bullet$, $\bullet \circ$, $\circ \bullet$, and $\circ \circ$
(the states of the PASEP on $2$ sites)
are $\frac{1}{q+5}$, $\frac{q+2}{q+5}$, $\frac{1}{q+5}$, and 
$\frac{1}{q+5}$, respectively; compare this with the right-hand
column of $SL_3$.  The steady-state probabilities of 
$\bullet \bullet \bullet$,
$\bullet \bullet \circ$,
$\bullet \circ \bullet$,
$\bullet \circ \circ$,
$\circ \bullet \bullet$,
$\circ \bullet \circ$,
$\circ \circ \bullet$,
$\circ \circ \circ$ are given by the right-hand column of $SL_4$.
This will be proved in Section \ref{PASEP}, building on work of 
\cite{CW}.

Note that since all coefficients of the matrix $SP_n$ are explicit (and
products of $q$-binomials) and since $SL_n$ comes from $SP_n$ by adding and
substracting rows, we have a simple expression of $E_I^J(q)$ as an alternating
sum of $C_I^J(q)$, hence of products of $q$-binomials.

Note that these matrices are invertible for generic values of $q$, in
particular for $q=0$. Hence, we can define Hall-Littlewood type functions by
\begin{equation}
\tilde{H}_J(q) = \sum_I E_I^J(q) L_I
\end{equation}
which interpolate between $S^I$ (at $q=1$) and a new kind of noncommutative
Schur functions $\Sigma_I$ at $q=0$.

We will see in the next section that the last column gives the $q$-enumeration
of permutation tableaux according to shape.
Let us write down the precise statement in that case.

\begin{proposition}
\label{ei1}
Let $c_I(q)$ be the coefficient of $\Psi_I$ in $S^{1^n}(q)$.
Let $e_I(q)$ be the coefficient of $L_I(q)$ in $S^{1^n}(q)$.

Then,
\begin{equation}
\label{eqei}
e_I(q) = \sum_{J\fatter I} (-1/q)^{l(I)-l(J)} q^{-\stb(I,J)} c_J(q), \text{ and }
\end{equation}
\begin{equation}
\label{eqci}
c_{j_1,\dots,j_r}(q) = 
[r]_q^{j_1} [r-1]_q^{j_2} \dots [2]_q^{j_{r-1}} [1]_q^{j_r}.
\end{equation}
\end{proposition}

\Proof
Straightforward from Theorem~\ref{recursive} and Equation~\ref{PsiL}.
\qed

We shall use the notation $\QStat_A(J):=c_J(q)$ in the sequel, regarding
it as a generalized $q$-factorial defined for all compositions (the classical
one coming from $J=(1^n)$.)

%%%%%%%%%%%%%%%%%%%%%%%%%%%%%%%%%%%%%
\subsubsection{A combinatorial lemma}

We need to describe the $q$-product in the $L(q)$ basis.
Our first objective will be to understand how $\st(I,K)$ can be related to
$\st(I,J)$ and $\st(J,K)$ for all $K\finer J\finer I$.
\begin{lemma}
\label{XYZ}
Let $X=\{x_1<\dots<x_r\} \subseteq Z=\{z_1<\dots<z_m\}$ be two sets of
positive integers.
For an integer $y$, let
\begin{equation}
\nu(y) = \#\{z\in Z|z\geq y\} + \#\{x\in X| x\leq y\},
\end{equation}
and for a set $Y$,
\begin{equation}
\nu(Y) = \sum_{y\in Y} \nu(y).
\end{equation}
For $r\leq s\leq m$, let
\begin{equation}
\Sigma_s(X,Z) =
\sum_{\genfrac{}{}{0pt}{}{X\subseteq Y\subseteq Z}{|Y|=s}}
      q^{\nu(Y)}.
\end{equation}
Then,
\begin{equation}
\Sigma_s(X,Z) = q^{\nu(X) + (r+1)(s-r)+\binom{s-r}{2}} \qbin{m-r}{s-r}.
\end{equation}
\end{lemma}

\Proof
Let $Z/X = U = \{ u_1<\dots<u_{m-r}\}$ and $\nu_i=\nu(u_i)$.
Then $\nu_{i}=m-i+1$, so that all $\nu_j$ are consecutive integers.
By definition,
\begin{equation}
\begin{split}
\Sigma_s(X,Z) &= \sum_{X\subseteq \{y_1<\dots<y_s\}\subseteq Z}
                q^{\nu(y_1)+\dots+\nu(y_s)}\\
&= q^{\nu(X)} \sum_{k_1<\dots<k_{s-r}} q^{\nu_{k_1}+\dots+\nu_{k_{s-r}}}\\
&= q^{\nu(X)} e_{s-r}(q^{\nu_1},\dots,q^{\nu_{m-r}}),
\end{split}
\end{equation}
where $e_{n}(X)$ is the usual elementary symmetric function on the alphabet
$X$. Thus,
\begin{equation}
\begin{split}
\Sigma_s(X,Z) &= q^{\nu(X)} e_{s-r}(q^{r+1},\dots,q^m) \\
&= q^{\nu(X)} q^{(r+1)(s-r)} e_{s-r}(1,q,\dots,q^{m-r-1}) \\
&= q^{\nu(X) + (r+1)(s-r)+\binom{s-r}{2}} \qbin{m-r}{s-r}.
\end{split}
\end{equation}
\qed

{\footnotesize
For example,
with $X=\{3,7\}$ and $Z=\{1,\dots,10\}$, one has
\begin{equation}
\Sigma_3(X,Z) = q^{\nu(3)+\nu(7)} \sum_{y\in Z/X} q^{\nu(y)}
= q^{15} (q^{10}+q^9+\dots+q^4+q^3) = q^{18} \qbin{8}{1}.
\end{equation}
\begin{equation}
\Sigma_4(X,Z) = q^{15} e_2(q^3,\dots,q^{10})
= q^{21} e_2(1,\dots,q^7) = q^{22} \qbin{8}{2}.
\end{equation}
}
%%%%%%%%%%%%%%%%%%%%%%%%%%%%%%%%%%%%%%%%%%%%%%%%%%%
\subsubsection{The $q$-product on the basis $L(q)$}

Recall that the $q$-product on $\NCSF$ is a non-associative product but that
$S^I(q)$ is  the  $q$-product on the parts of $I$, multiplied 
from right to left (see (\ref{non-assoc})).

\begin{lemma}
\label{LLP}
For an integer $p$ and a composition $I$, 
\begin{equation}
L_p(q) \qp L_I(q) = \sum_{K\finer p\mytr I}
   q^{st(K,p\mytr I)} \qbin{m_K+p}{p} \Psi_K,
\end{equation}
where $m_K = \#\{k\in\Des(K) | k\geq p\}$.
\end{lemma}

\Proof
Note that $L_p(q) = S_p$.  The $q$-products $S_p \qp \Psi_L$ are easily
computed by means of Lemma~\ref{qprod-init}.
\qed

We are now in a position to expand such $q$-products on the $L(q)$ basis:

\begin{lemma}
For an integer $p$ and a composition $I$, we have
\begin{equation}
\label{LLL}
L_p(q) \qp L_I(q) =
 \sum_{\genfrac{}{}{0pt}{}{J\finer p\mytr I}{j_1\geq p}}
     q^{st(J,p\mytr I) + \binom{l(J)-l(I)}{2} - \binom{l(I)}{2}}
     \qbin{l(I)+p-1}{l(J)-1} L_J(q).
\end{equation}
\end{lemma}

\Proof From Lemma \ref{LLP}, with $r=l(I)-1$, we have
\begin{equation}
\begin{split}
L_p(q) \qp L_I(q) &= \sum_{K\finer p\mytr I}
q^{st(K,p\mytr I)}
\sum_{s\ge r}q^{sr}\qbin{p-r}{p-s}\qbin{p+r}{s}\Psi_K\\
&=
\sum_{K\finer p\mytr I}\sum_{s\ge r}
\qbin{p+r}{s}
\left(
q^{st(K, p\mytr I)+rs}\qbin{p-r}{p-s}\Psi_K
\right)
\end{split}
\end{equation}
Thanks to Lemma~\ref{XYZ}, noting that
$\st(K,J)+\st(J,p\mytr I)=\nu(\Des(J))$
with $X=\Des(I)$ and $Z=\Des(K)\cap[p,\infty]$,
this is equal to
\begin{equation}
\begin{split}
&=
\sum_{s\ge r}\qbin{p+r}{s}
\sum_{\genfrac{}{}{0pt}{}{J\finer p\mytr I}{l(J)=s+1,\ j_1\ge p}}
q^{\st(J,p\mytr I)+\binom{s-r}{2}-\binom{r+1}{2}}
\sum_{K\succeq J}q^{\st(K,J)}\Psi_K\\
&=\sum_{\genfrac{}{}{0pt}{}{J\finer p\mytr I}{j_1\geq p}}
     q^{st(J,p\mytr I) + \binom{l(J)-l(I)}{2} - \binom{l(I)}{2}}
     \qbin{l(I)+p-1}{l(J)-1} L_J(q).
\end{split}
\end{equation}
\qed

{\footnotesize
For example,
\begin{equation}
L_{2}(q) \qp L_{21}(q) = [3] L_{41}(q) + [3] L_{311}(q)
 + [3] L_{221}(q) + q L_{2111}(q).
\end{equation}
\begin{equation}
L_{2}(q) \qp L_{12}(q) = [3] L_{32}(q) + q[3] L_{311}(q)
 + [3] L_{212}(q) + q^2 L_{2111}(q).
\end{equation}
\begin{equation}
\begin{split}
L_{3}(q) \qp L_{22}(q) =& [4] L_{52} + q\qbin{4}{2} L_{511}
                         +\qbin{4}{2} L_{412} + q^2[4] L_{4111}\\
&+ \qbin{4}{2} L_{322}(q) + q^2[4] L_{3211} + q[4] L_{3112}
 +q^4 L_{31111}.
\end{split}
\end{equation}
}

Since $S^J(q) = L_{j_1}(q) \qp (L_{j_2}(q)
                \qp (\dots (L_{j_{r-1}}(q) \qp L_{j_r}(q))\dots))$,
Formula~\ref{LLL} implies the following.

\begin{corollary}
The coefficient $E_I^J(q)$ is in $\NN[q]$.
\end{corollary}

\begin{corollary}
\label{rec-ei1}
Recall that $e_I(q)$ is the coefficient of $L_I(q)$ in $S^{1^n}(q)$.
Then, for any composition $I=(i_1,\dots,i_r)$,
\begin{equation}
\begin{split}
e_{1+i_1,i_2,\dots,i_r}(q)
& = [r]_q e_I + \sum_{k=1}^n q^{k-1} e_{i_1,\dots,i_k+i_{k+1},\dots,i_r}(q),\\
e_{1,i_1,i_2,\dots,i_r}(q)
& = e_{I}(q).
\end{split}
\end{equation}
Conversely, this property and the trivial initial conditions determine
completely the $e_{I}$.
\end{corollary}

\begin{proof}
This follows from the fact that $S^{1^n}(q) = 
L_{1}(q) \qp (\dots (L_{1}(q) \qp L_{1}(q))\dots)$, by 
putting $p=1$ into ~(\ref{LLL}).

\end{proof}

%%%%%%%%%%%%%%%%%%%%%%%%%%%%%%%%%%%%%%%%%%%%%%%%%%%%%%%%%%%%%%%%%%%%%
\subsubsection{Towards a combinatorial interpretation of  $E_I^J(q)$}

In Theorem~\ref{thm-Lig}, we will give a combinatorial interpretation
of the coefficients $E_I^J(q)$ expressing $S^J(q)$ in terms of the $L_I(q)$'s.
But first we need a new combinatorial algorithm sending a permutation to a
composition.

Let $\sigma$ be a permutation in $\SG_n$. We compute a composition 
$\Lig(\sigma)$ of $n$ as follows.

\begin{itemize}
\item Consider the Lehmer code of its inverse $\Lh(\sigma)$, that is, the word
whose $i$th letter is the number of letters of $\sigma$ to the left of $i$
and greater than $i$.

\item Fix $S=\emptyset$ and read $\Lh(\sigma)$ from right to left.
At each step, if the entry $k$ is strictly greater than the size of $S$, add
the ($k$-\#(S))-th element of the sequence $[1,n]$ with the
elements of $S$ removed .

\item The set $S$ is the descent set of a composition $C$, and
$\Lig(\sigma)$ is the mirror image $\bar C$ of  $C$.
\end{itemize}

{\footnotesize
For example, with $\sigma=(637124985)$, the Lehmer code of its
inverse is $\Lh(\sigma)=(331240010)$.
Then $S$ is $\emptyset$ at first, then the set $\{1\}$ (second step), then the
set $\{1,4\}$ (fifth step), then the set $\{1,4,2\}$ (eighth step).
Hence $C$ is $(1,1,2,5)$, so that $\Lig(\sigma)=(5,2,1,1)$.
}

One can find in Section~\ref{sec-patt} the permutations of $\SG_3$ and
$\SG_4$ arranged by rows according to their $\Lig$ statistics and by columns 
according to their recoil compositions.

%%%%%%%%%%%%%%%%%%%%%%%%%%%%%%%%%%%%%%
\subsubsection{A left $\Sym_q$-module}

Let $\sim$ be the equivalence relation on $\SG_n$ defined by
$\sigma\sim\tau$ whenever $\Lig(\sigma)=\Lig(\tau)$.
Let $\MM$ be the quotient of $\FQSym$~\cite{NCSF6} by the subspace
\begin{equation}
\VV=\{\F_\sigma -\F_\tau|\sigma\sim\tau\}\,.
\end{equation}
For a composition $I$, set $\kappa(I)=\maj(\bar I)$,
and for a permutation, $\kappa(\sigma)=\kappa(\Lig(\sigma))$.
Let $\LL_I$ denote the equivalence class of $q^{\kappa(\sigma)}\F_\sigma$.
Denote by $\circ_q$ be the $q$-product of $\FQSym$ inherited from $\WQSym$.
More precisely, if $\F'_\sigma=q^{\inv(\sigma)}\F_\sigma$ and
$\phi_q(\F_\sigma)=\F'_\sigma$, then
\begin{equation}
\F'_\sigma \, \circ_q \F'_\tau = \phi_q(\F_\sigma \F_\tau).
\end{equation}
This is the same structure as the one considered in \cite{NCSF6}.
In particular, in the basis $\G_\sigma=\F_{\sigma^{-1}}$, the product
is given by the $q$-convolution
\begin{equation}
\G_\alpha\circ_q\G_\beta=
\sum_{\gf{\gamma=u\cdots v}{\Std(u)=\alpha,\, \Std(v)=\beta}}
    q^{\inv(\gamma)-\inv(\alpha)-\inv(\beta)}\G_\gamma\,.
\end{equation}

\begin{lemma}
The quotient vector space $\MM$ is a left $\Sym$-module for the $q$-product of
$\FQSym$, that is, 
\begin{equation}
F\equiv G \mod \VV\ \Longrightarrow
 S_p\circ_q F\equiv S_p\circ_q G \mod \VV\,.
\end{equation}
\end{lemma}

\Proof
Let $\sigma^{-1}\sim\tau^{-1}\in\SG_l$ and $n=p+l$. We need to compare the
codes of the permutations appearing in the $q$-convolutions
\begin{equation}
U=\G_{12\cdots p} \circ_q \G_\sigma\
\text{and}\
V=\G_{12\cdots p}\circ_q \G_\tau\,.
\end{equation}
For a subset $S=\{s_1<s_2<\ldots <s_p\}$ of $[n]$, let $\sigma_S$
and $\tau_S$ be the elements of $U$ and $V$ whose prefix of length
$p$ is $s_1s_2\cdots s_p$. Then, the codes of $\sigma_S$ and
$\tau_S$ coincide on the first $p$ positions, and are equivalent
on the last $l$ ones, so that $\sigma_S\sim \tau_S$. Moreover,
$\sigma_S$ and $\tau_S$ arise with the same power of $q$, so
we have a module for the $q$-structure as well.
\qed

{\footnotesize
For example, 
\begin{equation}
\begin{split}
\F_{12} \circ_q q\F_{132}
&= q \F_{12354} + q^2 \F_{13254} + q^3 \F_{13524}
 + q^4 \F_{13542} + q^3 \F_{31254}\\
& + q^4 \F_{31524} + q^5 \F_{31542} + q^5 \F_{35124} + q^6 \F_{35142}
  + q^7 \F_{35412}\\
& = \LL_{41} + q\LL_{41} + \LL_{311}
 + \LL_{221} + q^2\LL_{41}\\
& + q\LL_{311}     + q\LL_{221}     + q^2 \LL_{311}  + q^2\LL_{221}
+ q \LL_{2111}.
\end{split}
\end{equation}
\begin{equation}
\begin{split}
\F_{12} \circ_q q \F_{312}
&= q \F_{12534} + q^2 \F_{15234} + q^3 \F_{15324}
 + q^4 \F_{15342} + q^3 \F_{51234}\\
& + q^4 \F_{51324} + q^5 \F_{51342} + q^5 \F_{53124} + q^6 \F_{53142}
  + q^7 \F_{53412}\\
& = \LL_{41} + q\LL_{41} + \LL_{311}
 + \LL_{221} + q^2\LL_{41}\\
& + q\LL_{311}     + q\LL_{221}     + q^2 \LL_{311}  + q^2\LL_{221}
+ q \LL_{2111}.
\end{split}
\end{equation}
\begin{equation}
\begin{split}
\F_{12} \circ_q q \F_{213}
&= q \F_{12435} + q^2 \F_{14235} + q^3 \F_{14325}
 + q^4 \F_{14352} + q^3 \F_{41235}\\
& + q^4 \F_{41325} + q^5 \F_{41352} + q^5 \F_{43125} + q^6 \F_{43152}
  + q^7 \F_{43512}\\
& = \LL_{41} + q\LL_{41} + \LL_{311}
 + \LL_{221} + q^2\LL_{41}\\
& + q\LL_{311}     + q\LL_{221}     + q^2 \LL_{311}  + q^2\LL_{221}
+ q \LL_{2111}.
\end{split}
\end{equation}
}

We have now:

\begin{lemma}
\label{lem-LL}
The left $q$-product of a $\LL_I$ by a complete function is given
by~(\ref{LLL}):
\begin{equation}
S_p\circ_q \LL_I = 
 \sum_{\genfrac{}{}{0pt}{}{J\finer p\mytr I}{j_1\geq p}}
     q^{st(J,p\mytr I) + \binom{l(J)-l(I)}{2} - \binom{l(I)}{2}}
     \qbin{l(I)+p-1}{l(J)-1} \LL_J.
\end{equation}
\end{lemma}

\Proof
Let us first show that this is true at $q=1$. Let $\sigma$ be such that
$\Lig(\sigma^{-1})=I$.
By definition of $\Lig$, the permutations $\tau$ occuring in
$\G_{12\cdots p}\circ_q\G_\sigma$ satisfy 
$\overline{\Lig(\tau^{-1})} \succeq \overline{\Lig(\sigma^{-1})}$, 
and the codes of those permutations have the form
\begin{equation}
s_1s_2\cdots s_p t_1t_2\cdots t_l\,,
\end{equation}
where $t=t_1t_2\cdots t_l$ is the code of $\sigma$
and $s_1\le s_2\le\ldots\le s_p$.
The compositions $J$ such that $l(J)-l(I)$ has a fixed value $m$ will
all be obtained by fixing the last $m$ values $s_p,\ldots,s_{p-m+1}$
in a way depending on the code $t$, the first $p-m$ being allowed to 
be any weakly increasing sequence
\begin{equation}
s_1\le s_2\le\ldots\le s_p\le l(J)-1\,,\ \text{which leaves }
\binom{p+l(I)-1}{l(J)-1}\ \text{choices.}
\end{equation}
Now, in the $q$-convolution $\G_{12\ldots p}\circ_q\G_\sigma$, these
permutations $\tau$ occur with a coefficient $q^{\inv(\tau)-\inv(\sigma)}$, so
that the coefficient of $\LL_J$ is, up to a power of $q$, the $q$-binomial
coefficient $\qbin{p+l(I)-1}{l(J)-1}$. By our choice of the
normalization $\LL_I=q^{\kappa(\sigma)}\F_\sigma$, this power of $q$ is the
same as in~(\ref{LLL}).
\qed

{\footnotesize
As one can check on the previous examples, we have indeed
\begin{equation}
\LL_{2} \circ_q \LL_{21} = (1+q+q^2) \LL_{41} + (1+q+q^2) \LL_{311}
 + (1+q+q^2) \LL_{221} + q \LL_{2111}.
\end{equation}
}

By Lemmas~\ref{lem-LL} and~(\ref{LLL}), the two bases $\LL$ and $L(q)$ have
the same multiplication formula, so that $E_I^J(q)$ is also the coefficient
of $\LL_I$ in the expansion of $S^J(q)$. Hence

\begin{theorem}
\label{thm-Lig}
Let $I$ and $J$ be two compositions of $n$.
Let $\PP(I,J)$ be the set of permutations
whose $\Lig$ statistic is $I$ and whose recoil composition is finer than $J$.
Then,
\begin{equation}
E_I^J(q) =
q^{-\maj(\overline{\Lig(\sigma)})} \sum_{\sigma\in \PP(I,J)} q^{\inv(\sigma)}.
\end{equation}
\end{theorem}

%%%%%%%%%%%%%%%%%%%%%%%%%%%%%%%%%%%%%%%%%%%%%%%%%%%%%%%%%%%%%%%%%%%%%%%%%%%%%%%
\subsection{The transition matrix $M(R(q),L(q))$}

The last transition matrix which remains to be computed is the one from $R(q)$
to $L(q)$.

%%%%%%%%%%%%%%%%%%%%%%%%%%%%%%
\subsubsection{First examples}

We have the following  matrices for $n=3,4$:

\begin{equation*}
RL_3 = M_3(R(q),L(q)) =
\left(
\begin{matrix}
1 & . & . & . \\
. & 1+q & 1 & . \\
. & . & 1 & . \\
. & . & . & 1 
\end{matrix}
\right)
\end{equation*}

\begin{equation*}
RL_4 =
\left(
\begin{matrix}
1 & . & . & . & . & . & . & . \\
. & 1+q+q^2 & 1+q & . & 1 & q & . & . \\
. & . & 1+q & . & 1 & . & . & . \\
. & . & q & 1+q+q^2 & . & 1+q & 1 & . \\
. & . & . & . & 1 & . & . & . \\
. & . & . & . & . & 1+q & 1 & . \\
. & . & . & . & . & . & 1 & . \\
. & . & . & . & . & . & . & 1
\end{matrix}
\right)
\end{equation*}

%%%%%%%%%%%%%%%%%%%%%%%%%%%%%%%%%%%%%%%%%%%%
\subsubsection{Combinatorial interpretation}

The coefficient of $L_I(q)$ in  $R_J(q)$  will be denoted
by $F_I^J(q)$.

>From the characterization  in Theorem~\ref{thm-Lig} of
$M(S(q),L(q))$ in terms of permutations 
we obtain:

\begin{theorem}
Let $I$ and $J$ be two compositions.
Let $\PP'(I,J)$ be the set of permutations whose $\Lig$ statistic is $I$ and
whose recoil composition is $J$.
The coefficient $F_I^J$ of $L_I(q)$ in the expansion of $R_J(q)$
is given by
\begin{equation}
q^{-\maj(\overline{\Lig(\sigma)})} \sum_{\sigma\in \PP'(I,J)} q^{\inv(\sigma)}.
\end{equation}
\end{theorem}

%%%%%%%%%%%%%%%%%%%%%%%%%%%%%%%%%%%%%%%%%%%%%%%%%%%%%%%%%%%%%%%%%%%%%%%%%%%%%%%
%%%%%%%%%%%%%%%%%%%%%%%%%%%%%%%%%%%%%%%%%%%%%%%%%%%%%%%%%%%%%%%%%%%%%%%%%%%%%%%
%%%%%%%%%%%%%%%%%%%%%%%%%%%%%%%%%%%%%%%%%%%%%%%%%%%%%%%%%%%%%%%%%%%%%%%%%%%%%%%
\section{The PASEP and type A permutation tableaux}\label{PASEP}

Permutation tableaux (of type A) are certain fillings of Young diagrams with
$0$'s and $1$'s which are in bijection with permutations (see \cite{SW} for
two bijections).
They are a distinguished subset of Postnikov's (type A)
$\Le$-diagrams~\cite{Postnikov}, which index cells of the totally non-negative
part of the Grassmannian.

Apart from this geometric connection, permutation tableaux are of interest as
they are closely connected to a model from statistical physics called the
partially asymmetric exclusion process (PASEP) \cite{CW}.
More precisely, the PASEP with $n$ sites is a model in which particles hop
back and forth (and in and out) of a one-dimensional lattice, such that at
most one particle may occupy a given site (the probability of hopping left
is $q$ times the probability of hopping right.)
See \cite{CW} for full details. Therefore there are $2^n$ possible states
of the PASEP. There is a simple bijection from a state $\tau$ of the PASEP to
a Young diagram $\lambda(\tau)$ whose semiperimeter is $n+1$.
The main result of \cite{CW} is that the steady state probability that the
PASEP is in configuration $\tau$ is equal to the $q$-enumeration of
permutation tableaux of shape $\lambda(\tau)$ divided by the $q$-enumeration
of all permutation tableaux of semiperimeter $n+1$.

In this section we will give an explicit formula for the $q$-enumeration of
permutation tableaux of a given shape.
So in particular this is an explicit formula for the steady state probability
of each state of the PASEP. Additionally, by results of~\cite{SW}, this
formula counts permutations with a given set of weak excedances according to
{\it crossings}; it also counts permutations with a given set of {\it descent
bottoms} according to occurrences of the pattern $2-31$.

%%%%%%%%%%%%%%%%%%%%%%%%%%%%%%%%%%%%%%%%%%%%%%%%%%%%%%%%%%%%%%%%%%%%%%%%%%%%%%%
\subsection{Permutation tableaux}

Regard the following $(k, n-k)$ rectangle (here $k=3$ and $n=8$)
\begin{equation}
\tableaux{\\{}&{}&{}&{}&{} \\
{}&{}&{}&{}&{} \\
{}&{}&{}&{}&{}\\&}
\end{equation}
as a poset $Q^A_{k,n}$: the elements of the poset are the boxes, and box $b$
is less than $b'$ if $b$ is southwest of $b'$.  We then define a
{\it type $A$ Young diagram} contained in a $(k,n-k)$ rectangle
to be an order ideal in the poset $Q^A_{k,n}$.
This corresponds to the French notation for representing Young diagrams.
We will sometimes refer to such a Young diagram by the partition $\lambda$
given by the lengths of the rows of the order ideal.  Note that we allow
partitions to have parts of size $0$.

As in \cite{SW}, we define a type A {\em permutation tableau} $\T$ to be a
type A Young diagram $Y_\lambda$ together with a filling of the boxes with
$0$'s and $1$'s such that the following properties hold:
\begin{enumerate}
\item Each column of the diagram contains at least one $1$.
\item There is no $0$ which has a $1$ below it in the same column
{\em and} a $1$ to its left in the same row.
\end{enumerate}

\noindent
We call such a filling a \emph{valid} filling of $Y_\lambda$.
Here is an example of a type A permutation tableau.
\begin{equation}
\tableaux{\\{0}&{1}&{1}&& \\
{1}&{0}&{1}&{1}& \\
{1}&{0}&{1}&{0}&{1}\\&}
\end{equation}

\noindent
Note that if we forget the requirement (1) in the definition of type A
permutation tableaux then we recover the description of a (type A)
$\Le$-diagram~\cite{Postnikov}, an object which represents a cell in the
totally nonnegative part of a Grassmannian. In that case, the total number of
$1$'s corresponds to the dimension of the cell.

We define the {\it rank} $\rk(\T)$ of a permutation tableau (of type A) $\T$
with $k$ columns to be the total number of $1$'s in the filling minus $k$.
(We subtract $k$ since there must be at least $k$ $1$'s in a valid filling of
a tableau with $k$ columns.)

%%%%%%%%%%%%%%%%%%%%%%%%%%%%%%%%%%%%%%%%%%%%%%%%%%%%%%%%%%%%%%%%%%%%%%%%%%%%%%%
\subsection{Enumeration of permutation tableaux by shape}

%%%%%%%%%%%%%%%%%%%%%%%%%%%%%%%%%%%%%%%%%%%%%%%%%%%%%%%%%%%%%%%%%%%%%%%%%%%%

Starting from a partition with $k$ rows and $n-k$ columns,
one encodes it as a composition $I=(i_1,\dots,i_k)$ of $n$ as follows:
$i_1-1$ is the number of columns of length $k$,
$i_2-1$ is the number of columns of length $k-1$, and
\dots, $i_k-1$ is the number of columns of length $1$.

Let $\ell(I)$ denote the number of parts of $I$.
Then the number $\PT^A_I$ of permutation tableaux of shape corresponding to
$I$ is given by a simple formula coming from combinatorics of noncommutative
symmetric functions.
Indeed, according to~\cite[Proposition 9.2]{Tev},
\begin{equation}
L_1^n = \sum_{I\vDash n} g_I \Psi_I,
\end{equation}
where
\begin{equation}
g_I = \prod_{k=1}^{l(I)} (l(I)-k+1)^{i_k}.
\end{equation}
Hence, the coefficient $e_J$ of
\begin{equation}
L_1^n = \sum_{J\vDash n} e_J L_J
\end{equation}
is given by
\begin{equation}
e_J = \sum_{I\finer J} (-1)^{l(I)-l(J)} \prod_{k=1}^{l(I)} (l(I)-k+1)^{i_k}.
\end{equation}
Moreover, from~\cite{HNTT}, Theorem 5.1, we known that $e_I$ is the number of
permutations such that $\GC(\sigma)=I$. Finally, since permutation tableaux of
a given shape are in bijection with permutations with given descent
bottoms~\cite{SW} and that $\GC$ does the same up to reverse complement of the
permutations, this number is also the number of permutation tableaux of shape
$I$.

\begin{theorem}
\label{Theorem1}
\begin{equation}
\PT^A_I = \sum_{J\fatter I} (-1)^{\ell(I)-\ell(J)} \Stat(J),
\end{equation}
where the sum is over the compositions $J$ coarser than $I$ and where
$\Stat$ is defined by
\begin{equation}
\Stat(j_1,\dots,j_p) := p^{j_1} (p-1)^{j_2} \dots 2^{j_{p-1}} 1^{j_p}.
\end{equation}
\end{theorem}
\qed

For example with $I=(3,4,1)$,
we get
\begin{equation}
\PT^A_{341} = 3^3 2^4 1^1 - 2^7 1^1 - 2^3 1^5 + 1^8 = 297.
\end{equation}

%%%%%%%%%%%%%%%%%%%%%%%%%%%%%%%%%%%%%%%%%%%%%%%%%%%%%%%%%%%%%%%%%%%%%%%%%%%%
\subsubsection{$q$-enumeration of permutation tableaux according to their
shape}

In this section, we make the connection between the coefficients $e_I(q)$
previously seen, and the $q$-enumeration of permutation tableaux.
Recall that $e_I(q)$ is the coefficient of $L_I(q)$ in $S^{1^n}(q)$.
We saw in Corollary \ref{rec-ei1} that for all compositions $I=(i_1,\dots,i_r)$,
the following hold:
\begin{itemize}
\item $e_{(1,i_1,i_2,\dots,i_r)}(q)
 = e_{I}(q).$
\item $e_{(1+i_1,i_2,\dots,i_r)}(q)
 = [r]_q e_I + \sum_{k=1}^{r-1} q^{k-1}
e_{(i_1,\dots,i_k+i_{k+1},\dots,i_r)}(q)$
\end{itemize}

\noindent
It is possible to transform this result into a $q$-enumeration of
permutation tableaux by their rank.
Let
\begin{equation}
\PT^A_I(q) := \sum_T q^{\rk(T)},
\end{equation}
where the sum is over all permutation tableaux whose shape corresponds to $I$.

The following result generalizes Theorem \ref{Theorem1}.
Its proof follows directly from Proposition~\ref{ei1},
Corollary~\ref{rec-ei1}, and Lemma~\ref{lem-ptA} below.

\begin{theorem}
\label{Theorem2}
Let $I$ be a composition. Then,
\begin{equation}
\label{eqptA}
\PT^A_I(q)
= e_I(q)
= \sum_{J\fatter I} (-1/q)^{l(I)-l(J)} q^{-\stb(I,J)} \QStat_A(J),
\end{equation}
where $\QStat_A$ is recalled to be
\begin{equation}
\QStat_A(j_1,\dots,j_p) := [p]_q^{j_1} [p-1]_q^{j_2} \dots [2]_q^{j_{p-1}}
[1]_q^{j_p}.
\end{equation}
\end{theorem}

By the results of~\cite{CW}, Theorem \ref{Theorem2} gives an explicit 
formula for the steady
state probabilities in the partially asymmetric exclusion process (PASEP).
More specifically, consider the PASEP on a one-dimensional lattice of $n$
sites where particles hop right with probability $dt$, hop left with
probability $q dt$, enter from the left at a rate $dt$, and exit to the right
at a rate $dt$. Let us number the $n$ sites from {\it right to left} with the
numbers $1$ through $n$.  Then we have the following result.

\begin{corollary}
Recall the notation of Theorem~\ref{Theorem2}.  Let $I$ be a composition of
$n+1$, and let $Z_n$ denote the partition function for the PASEP.
Let $\tau$ denote the state of the PASEP
in which all sites of $\Des(I)$ are occupied
by a particle
and all sites of $[n-1]\setminus \Des(I)$ are empty.
Then the probability that in the steady state, the PASEP is in 
state $\tau$, is 
\begin{equation*}\frac{\sum_{J\fatter I} (-1/q)^{l(I)-l(J)} q^{-\stb(I,J)} 
\QStat_A(J)}{Z_n}.
\end{equation*}
\end{corollary}

By the results of \cite{SW}, this is also  an explicit formula 
enumerating permutations with a fixed set of {\it weak excedances}
according to the number of {\it crossings}; equivalently,
an explicit formula enumerating permutations with a fixed
set of {\it descent bottoms} according to the number of occurences
of the {\it generalized pattern} $2-31$.  See \cite{SW} 
for definitions.

More specifically,
let $I$ be a composition of $n+1$, let $DB(I)$ be the 
descent set of the reverse composition of $I$, 
and let $W(I) = \{1\} \cup \{1 + DB(I) \}$.  Here
$1+DB(I)$ denotes the set obtained by adding $1$ to each element
of $DB(I)$.   If $\sigma$ is a permutation,
let $(2-31)\sigma$ denote the number of occurrences of the pattern
$2-31$ in $\sigma$, and let $\Cr(\sigma)$ denote the number of 
{\it crossings} of $\sigma$.
Let $T_I(q) = \sum_{\sigma} q^{(2-31)\sigma}$ be the 
sum over all permutations in $S_{n+1}$ whose set of descent 
bottoms is $DB(I)$.  And let $T'_I(q) = \sum_{\sigma} q^{\Cr(\sigma)}$
be the sum over all permutations in $S_{n+1}$ whose 
set of weak excedances is $W(I)$.

\begin{corollary}
\begin{equation*}T_I(q) = T'_I(q) = 
{\sum_{J\fatter I} (-1/q)^{l(I)-l(J)} q^{-\stb(I,J)} 
\QStat_A(J)}.
\end{equation*}

\end{corollary}

\medskip
For example with $I=(3,4,1)$, the compositions coarser than $I$ are
$(3,4,1)$, $(7,1)$, $(3,5)$, and $(8)$, so
we get
\begin{equation}
\begin{split}
\PT^A_{341}(q)
&= \frac{1}{q^2}\left(  \frac{[3]_q^3 [2]_q^4}{q}
                                - \frac{[2]_q^7}{q}
                                - \frac{[2]_q^3}{1}
                                + \frac{1}{1}
\right)\\
&=
q^7  + 7 q^6  + 24 q^5  + 52 q^4  + 76 q^3  + 75 q^2  + 47 q + 15.
\end{split}
\end{equation}

The descent set $D(I)$ of $I$ is $\{3,7\}$, which corresponds to 
the following state of the PASEP:
$\tau=\bullet \circ \circ \circ \bullet \circ \circ$.  
Therefore the probability that in the steady state, the PASEP
is in state $\tau$, is 
$\frac{q^7  + 7 q^6  + 24 q^5  + 52 q^4  + 76 q^3  +75q^2 + 47 q + 15}{Z_7}$.

The polynomial 
$q^7  + 7 q^6  + 24 q^5  + 52 q^4  + 76 q^3  +75q^2 + 47 q + 15$
also enumerates the permutations in $S_8$ 
with set of descent bottoms $\{1,5\}$ according to 
occurrences of the pattern $2-31$.  And it enumerates permutations
in $S_8$  
with weak excedances in positions $\{1,2,6\}$ according to crossings.

\medskip

The reader might want to compare~(\ref{eqptA}) with~(\ref{eqei}).

\begin{lemma}
\label{lem-ptA}
Let $I=(i_1,\dots,i_r)$ be a composition. Then
\begin{equation}
\PT^A_{(1,i_1,i_2,\dots,i_r)}(q) = \PT^A_{I}(q),
\end{equation}
\begin{equation}
\PT^A_{(1+i_1,i_2,\dots,i_r)}(q)
   = [r]_q \PT^A_I
     + \sum_{k=1}^n q^{k-1} \PT^A_{(i_1,\dots,i_k+i_{k+1},\dots,i_r)}(q).
\end{equation}
\end{lemma}

\Proof
First note that
$PT^A_{(1,i_1,i_2,\dots,i_r)}(q)= PT^A_{I}(q)$: this just says that the
$q$-enum\-eration of permutation tableaux of shape $\lambda$ is the same as
the $q$-enumeration of permutation tableaux of shape $\lambda'$, where
$\lambda'$ is obtained from $\lambda$ by adding a row of length $0$.

Therefore we just need to prove the second equality.
Let $\lambda=(\lambda_1,\dots,\lambda_r)$ be a partition. Then, in terms of
partitions, the statement translates as:
\begin{equation}
\PT^A_{\lambda}(q)
 = [r] PT^A_{(\lambda_1-1, \lambda_2-1, \dots, \lambda_r-1)}(q) + 
   \sum_{k=1}^r q^{k-1} PT^A_{(\lambda_1,\dots,\lambda_{r-k},
      \widehat{\lambda_{r-k+1}},\lambda_{r-k+2}-1,\dots,\lambda_{r}-1)}(q),
\end{equation}
where $\widehat{\lambda_{r-k+1}}$ means that this part has been removed.

To this aim, we need to introduce the notion of a {\it restricted} zero.  We
say that a zero in a tableau is {\it restricted} if there is a $1$ below it
in the same column. Note that every entry to the left of and in the same row
as the restricted zero must also be zero.

We will prove the recurrence by examining the various possibilities for the
set $S$ of $r$ boxes of the Young diagram $\lambda$ which are rightmost in
their row. We will partition (most of) the permutation tableaux with shape
$\lambda$ based on the position of the highest restricted zero among $S$.

\begin{figure}[ht]
\centering
\includegraphics[height=1.8in]{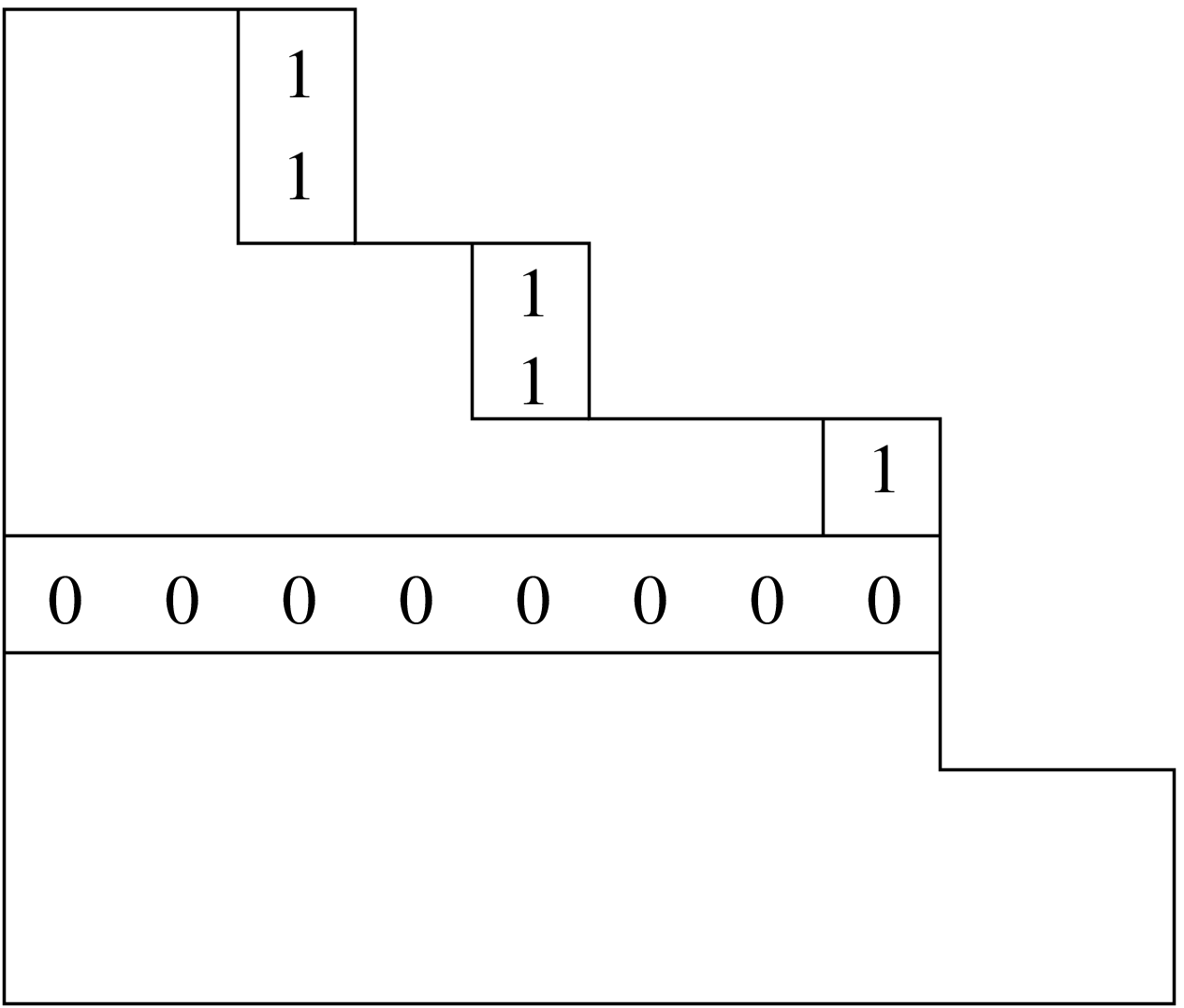}
\caption{}
\label{Step1}
\end{figure}

We will label rows of the Young diagram from top to bottom, from $1$ to $r$.
Consider the set of tableaux obtained via the following procedure: choose a
row $k$ for $1 \leq k \leq r-1$, and fill it entirely with $0$'s.  Also fill
each box of $S$ in row $\ell$ for any $\ell>k$ with a $1$.  Now ignore row
$k$ and the filled boxes of $S$, and fill the remaining boxes (which can be
thought of as boxes of a partition of shape
$\lambda':=(\lambda_1,\dots,\lambda_{k-1},\widehat{\lambda_k},
\lambda_{k+1}-1,\dots,\lambda_r-1)$) in any way which gives a legitimate
permutation tableau of shape $\lambda'$ (see Figure \ref{Step1}.)
Note that if we add back the ignored boxes, we will increase the rank of the
first tableau by $k-1$. So the $q$-enumeration of the tableaux under
consideration is exactly
$\sum_{k=1}^r q^{k-1} PT_{(\lambda_1,\dots,\lambda_{r-k},
\widehat{\lambda_{r-k+1}},\lambda_{r-k+2}-1,\dots,\lambda_{r}-1)}(q).$

Let us denote the columns of $\lambda$ which contain a north-east corner of
the Young diagram $\lambda$ as $c_1,\dots,c_h$; we will call them
\emph{corner columns}. Denote the lengths of those columns by $C_1,\dots,C_h$,
so $C_1>\dots >C_h$.  And denote the differences of their lengths by
$d_1:=C_1-C_2,\dots,d_{h-1}:=C_{h-1}-C_h,d_h:=C_h$.

Clearly, our procedure constructs all permutation
tableaux of shape $\lambda$ with the following description: at least one box
of $S$ is a restricted zero.  Furthermore, if we choose the restricted zero of
$S$ (say in box $b$) which is in the lowest row (say row $k$), then every box
of $S$ in a row above $k$ is filled with a $1$. Equivalently, each corner
column $c_j$ left of $b$ has its top $d_j$ boxes filled with $1$'s, and
contains at least $d_j+1$ ones total; and the corner column containing $b$
contains at least $d+1$ ones total, where $d$ is the number of boxes above
$b$ in the same column.

The permutation tableaux of shape $\lambda$ which this procedure has {\it not
constructed} are those tableaux such that either no box of $S$ is a restricted
zero, or else there {\it is} a box of $S$ which is a restricted zero.
Let $b$ denote the lowest such box. The condition that all boxes of $S$ above
$b$ must be $1$'s is violated.
Let $W$ denote this set of tableaux.

\begin{figure}[ht]
\centering
\includegraphics[height=1.8in]{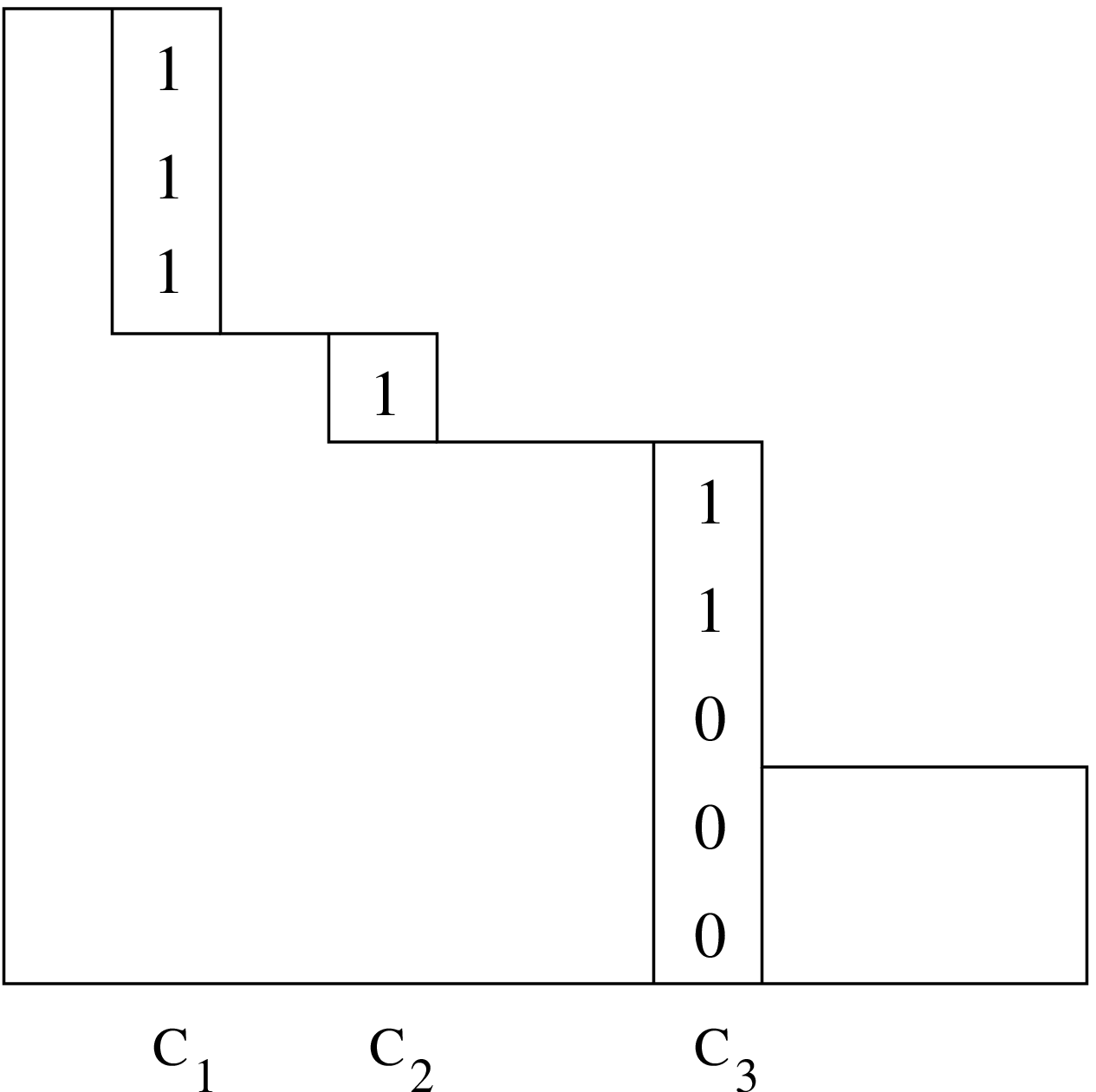}
\caption{}
\label{Step2}
\end{figure}

The following construction gives rise to all permutation tableaux in $W$ (See
Figure \ref{Step2}.)
Choose a corner column $c_j$ and a number $m$ such that $1 \leq m \leq d_j$.
Fill the top $m$ boxes of $c_j$ with $1$'s and the remaining boxes with $0$'s.
For each $i<j$, fill the top $d_i$ boxes of column $c_i$ with $1$'s.  Now
ignore the boxes that have been filled, and choose any filling of the
remaining boxes -- which form a partition of shape
$\lambda'':=(\lambda_1-1,\lambda_2-1,\dots,\lambda_r-1)$
-- 
which gives a legitimate permutation tableau of shape $\lambda''$.
Note that adding back the boxes we had ignored will add $d_1+\dots+d_{h-1}+m$
to the  rank of the tableau of shape $\lambda'$. 
Since the quantity $d_1+\dots+d_{h-1}+m$ can range between $0$ and $r-1$,
the rank of the tableaux in $W$ is
$ [r] PT_{(\lambda_1-1, \lambda_2-1, \dots, \lambda_r-1)}(q)$. 
\qed

\begin{figure}[ht]
\centering
\includegraphics[height=1.8in]{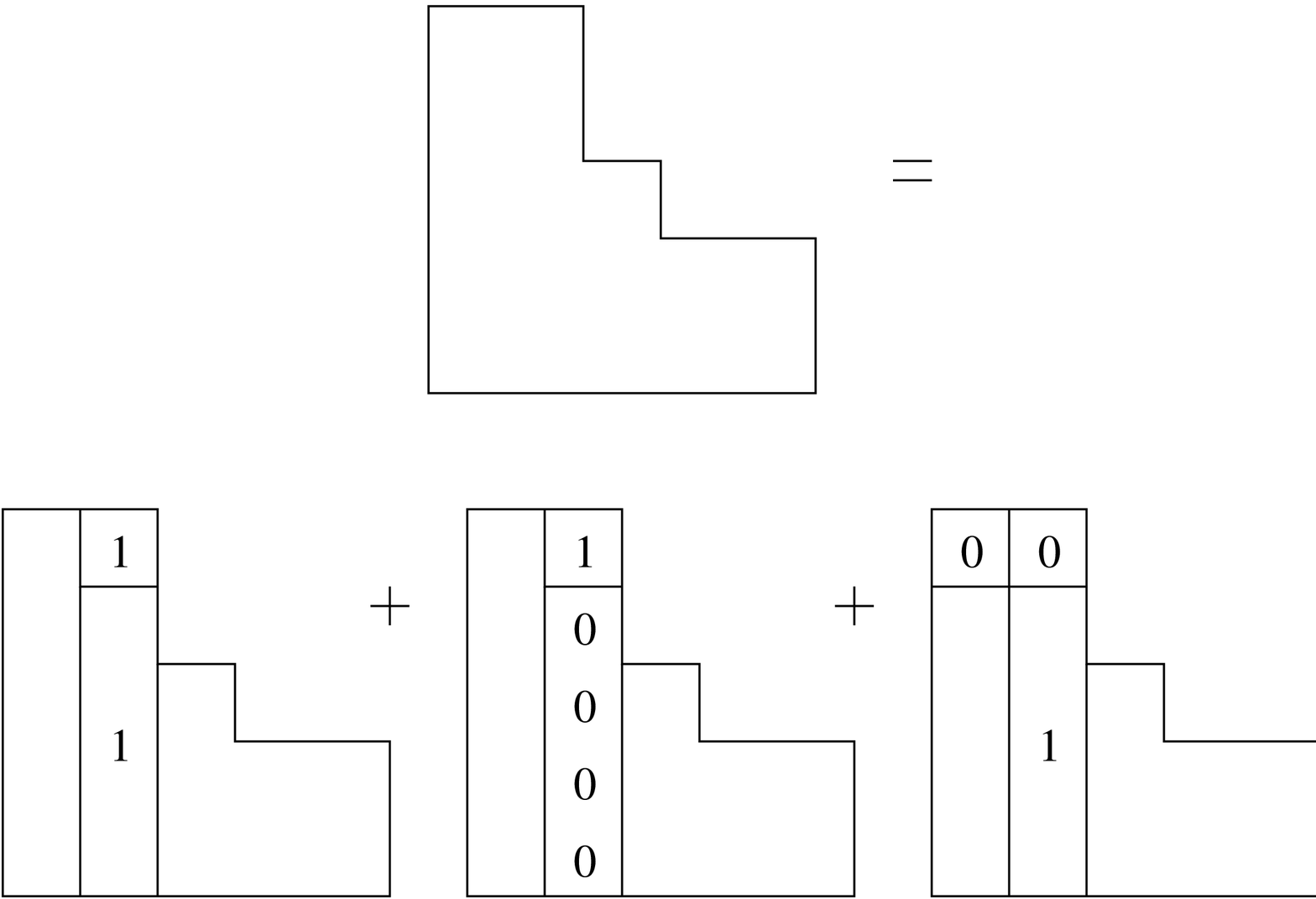}
\caption{}
\label{ARecur}
\end{figure}

Note that we could give an alternative (direct) proof of
Theorem~\ref{Theorem2} by using the following recurrences for permutation
tableaux (which had been observed in~\cite{SW}).
See Figure \ref{ARecur} for an illustration of the second recurrence.

\begin{lemma}\label{useful}
The following recurrences for type A permutation tableaux hold.\\
\begin{itemize}
\item $PT^A_{(i_2,i_3,\dots,i_n)}(q) = PT^A_{(1,i_2,i_3,\dots,i_n)}(q)$\\
\item
     $PT^A_{(i_1,i_2,\dots,i_n)}(q) = q PT^A_{(i_1 - 1,
i_2+1,i_3,\dots,i_n)}(q) +
   PT^A_{(i_1 -1,i_2,\dots,i_n)}(q)\\ +
PT^A_{(1,i_1+i_2-1,i_3,\dots,i_n)}(q)$.
\end{itemize}
\end{lemma}

%%%%%%%%%%%%%%%%%%%%%%%%%%%%%%%%%%%%%%%%%%%%%%%%%%%%%%%%%%%%%%%%%%%%%%%%%%%%%%%
%%%%%%%%%%%%%%%%%%%%%%%%%%%%%%%%%%%%%%%%%%%%%%%%%%%%%%%%%%%%%%%%%%%%%%%%%%%%%%%
%%%%%%%%%%%%%%%%%%%%%%%%%%%%%%%%%%%%%%%%%%%%%%%%%%%%%%%%%%%%%%%%%%%%%%%%%%%%%%%
\section{Permutation tableaux and enumeration formulas in type B}

One can also define~\cite{LW} type B $\Le$-diagrams and permutation tableaux,
where the Type B $\Le$-diagrams index cells in the odd orthogonal
Grassmannian, and type B permutation tableaux are in bijection with signed
permutations.
In this section, we will enumerate permutation tableaux of type B of a fixed
shape, according to rank.  This formula can be given an interpretation in
terms of signed permutations.

To define type $B_n$ Young diagrams, regard the following shape
\begin{equation}
\tableaux{\\{}&{}&{}&{}\\{}&{}&{}\\{}&{}\\{}\\&}
\end{equation}
as representing a poset $Q^B_n$ (here $n=4$):
the elements of the poset are the boxes, and box $b$ is less than
$b'$ if $b$ is southwest of $b'$.
We then define a {\it type $B_n$ Young diagram} to be an order ideal in the
poset $Q^B_n$.

As in \cite{LW}, we define a type B {\em permutation tableau} $\T$ to be a
type B Young diagram $Y_\lambda$ together with a filling of the boxes with
$0$'s and $1$'s such that the following properties hold:
\begin{enumerate}
\item Each column of the diagram contains at least one $1$.
\item There is no $0$ which has a $1$ below it in the same column
{\em and} a $1$ to its left in the same row.
\item If a diagonal box contains a $0$,  every box in that
 row must contain a $0$.
\end{enumerate}

\noindent
Here is an example of a type B permutation tableau.
\begin{equation}
\tableaux{\\{1}\\{0}&{0}&{1}\\
{0}&{0}&{0}\\
{1}&{1}\\
{1}\\&}
\end{equation}

Note that if we forget requirement (1) in the definition of a type B
permutation tableaux then we recover the description of a type B
$\Le$-diagram~\cite{LW}, an object which represents a cell in the totally
nonnegative part of an odd orthogonal Grassmannian.

As before, we define the \emph{rank} $\rk(\T)$ of a permutation tableau $\T$
(of type B) with $k$ columns to be the total number of $1$'s in the filling
minus $k$.

Starting from a type B Young diagram $Y_{\lambda}$ inside a staircase of
height $n+1$, we encode it as a composition of $n$ as follows.
If $k$ is the width of the widest row of $Y_{\lambda}$, then
$I=(i_1,\dots,i_{k+1})$ is defined by:
$i_1+1$ is the number of rows of length $k$,
$i_2$ is the number of rows of length $k-1$,
\dots, $i_{k+1}$ is the number of rows of length $0$.

We now explain how to enumerate type B permutation tableaux
of a fixed shape according to their rank.

Define $\QStat_B$ by
\begin{equation}
\begin{split}
\QStat_B(j_1,\dots,j_p)
 :=& \QStat_A(j_1,\dots,j_p) \prod_{t=1}^{p-1} \left(1+q^{t}\right) \\
=& [p]_q^{j_1} [p-1]_q^{j_2} \dots [2]_q^{j_{p-1}} [1]_q^{j_p}
  \prod_{t=1}^{p-1} \left(1+q^{t}\right).
\end{split}
\end{equation}

\begin{theorem}
\label{TheoremB}
Let $I$ be a composition.
\begin{equation}
\label{eqptB}
\PT^B_I(q)
= \sum_{J\fatter I} (-1/q)^{l(I)-l(J)} q^{-\stb(I,J)} \QStat_B(J).
\end{equation}
where $p$ is the length of $J$.
\end{theorem}

Note that the formula enumerating type B permutation tableaux is
very similar to the formula enumerating type A permutation tableaux.

As an example, suppose we want to enumerate according to rank the type $B$
permutation tableaux that have the following shape:
\begin{equation}
\tableaux{\\{}\\{}&{}\\
{}&{}\\
{}\\
&}
\end{equation}
We take $n=4$ (we would get the same answer for any $n>4$), and $k=2$ since
the widest row has width $2$.  Then the corresponding composition is
$I=(1,2,0)$.
We then get

\begin{align*}
\PT^B_{(1,2,0)}(q)&=q^{-2}(q^{-1}[3]_q [2]_q^2  (1+q) (1+q^2)
-q^{-1}[2]_q^{3}(1+q) -
[2]_q  (1+q)
+1)\\ &= q^4+4q^3+8q^2+10q+6.
\end{align*}

We will prove Theorem \ref{TheoremB} directly:
we first prove some recurrences
for type B permutation tableaux, and then prove that the formula in
Theorem \ref{TheoremB} satisfies the same recurrences.
\begin{lemma}\label{Brecur}
The following recurrences for type B permutation tableaux hold.
\begin{equation}
\PT^B_{(0,i_2,i_3,\dots,i_k)}(q) = \PT^B_{(i_2,i_3,\dots,i_k)}(q),\\
\end{equation}
\begin{equation}
\begin{split}
\PT^B_{(i_1,i_2,i_3,\dots,i_k)}(q) =&
\PT^B_{(i_1-1,i_2,i_3,\dots,i_k)}(q)+ qPT^B_{(i_1-1,i_2+1,i_3,\dots,i_k)}(q)\\
&+ PT^B_{(0,i_1+i_2-1,i_3,\dots,i_k)}(q).
\end{split}
\end{equation}
\end{lemma}

\Proof
The first recurrence says that enumerating permutation tableaux of a shape
which has a unique row of maximal width is the same as enumerating permutation
tableaux of the shape obtained from the first shape by deleting the rightmost
column.  This is clear, since the rightmost column will have only one box
which must be filled with a $1$.

\begin{figure}[ht]
\centering
\includegraphics[height=1.5in]{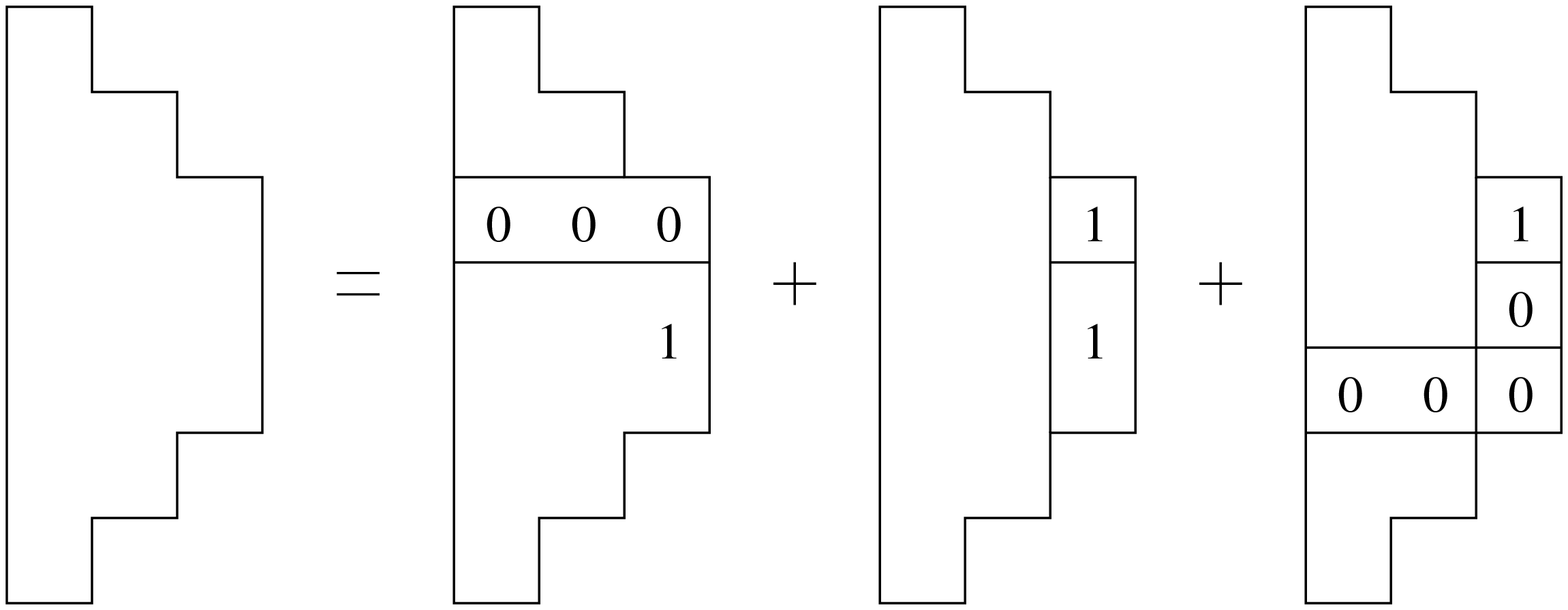}
\caption{}
\label{BRecur}
\end{figure}

To see that the second recurrence holds, see Figure \ref{BRecur}.
Consider the topmost box $b$ of the rightmost column of an arbitrary type B
permutation tableau of shape corresponding to $(i_1,\dots,i_k)$. Since
$i_1>0$, the rightmost column has at least two boxes.
If $b$ contains a $0$, then by definition of type B permutation tableaux,
there is a $1$ below it in the same column -- which implies that the entire
row containing $b$ must be filled with $0$'s.  We can delete that entire row
and what remains will be a type B permutation tableau (of smaller shape).

If $b$ contains a $1$ and there is another $1$ in the same column, then we can
delete the box $b$ and what remains will be a type B permutation tableau.

If $b$ contains a $1$ and there is no other $1$ in the same column, then let
$b'$ denote the bottom box of that column.  By the definition of type B
permutation tableau, the entire row of $b'$ is filled with $0$'s. If we
delete the entire row of $b'$ and every box below and in the same column as
$b$, then what remains will be a type B permutation tableau.
\qed

Now we prove Theorem \ref{TheoremB}.

\Proof[of the theorem]

Let
\begin{equation}
g_I(q) := 
\sum_{J\fatter I} (-1/q)^{\ell(I)-\ell(J)} q^{-\stb(I,J)} \QStat_B(J).
\end{equation}
We want to prove that $\PT^B_I(q) = g_I(q)$.
We claim that it is enough to prove the following two facts:
\begin{enumerate}
\item $g_{(0,i_2,i_3,\dots,i_n)}(q) = g_{(i_2,i_3,\dots,i_n)}(q)$
\item $g_{(i_1,i_2,\dots,i_n)}(q) = q g_{(i_1-1,i_2+1,i_3,\dots,i_n}(q) +
          g_{(i_1 -1,i_2,i_3,\dots,i_n)}(q) +
g_{(0,i_1+i_2-1,i_3,\dots,i_n)}(q)$
       when $i_1>0$.
\end{enumerate}
By Lemma \ref{Brecur}, both of these recurrences are true for $PT^B_I(q)$.
And the two recurrences together clearly determine $g_{I}(q)$ for any
composition $I$, which is why it suffices to prove these recurrences.

Consider the first recurrence.
To prove it, we will pair up the terms that occur in
\begin{equation}
g_{0,i_2,\dots,i_n}(q) := \sum_{J\fatter I} (-1/q)^{\ell(I)-\ell(J)}
q^{-\stb(I,J)} \QStat_B(J),
\end{equation}
pairing each composition of the form
$J:=(0,j_1,j_2,\dots,j_r)$ with the composition $J':=(j_1,j_2,\dots,j_r)$.

Note that $\ell(J)=\ell(J')+1$ and $\stb(I,J) = \stb(I,J')+1$.
Also
\begin{equation}
\QStat_B(0,j_1,j_2,\dots,j_r) = \QStat_B(j_1,j_2,\dots,j_r) (1+q^r)
\end{equation}
so that
\begin{equation}
\QStat_B(0,j_1,j_2,\dots,j_r) - \QStat_B(j_1,j_2,\dots,j_r) =
q^{r} \QStat_B(j_1,j_2,\dots,j_r).
\end{equation}

And now it follows from the fact that
\begin{equation}
\stb((0,I),(0,J)) = q^{r-1} \stb(I,J),
\end{equation}
that the contribution to $g_{0,i_2,\dots,i_n}(q)$ by the pair of compositions
$J$ and $J'$ is exactly the contribution to $g_{i_2,\dots,i_n}(q)$ by the
composition $(j_1,\dots,j_r)$.
So $g_{(0,i_2,i_3,\dots,i_n)}(q) = g_{(i_2,i_3,\dots,i_n)}(q)$.

\bigskip
Now let us turn our attention to the second recurrence.
We prove the second recurrence by showing that each term of
$g_{(i_1,\dots,i_n)}(q)$ comes from either one term each from
$qg_{(i_1-1,i_2+1,i_3,\dots,i_n)}(q)$ and
$g_{(i_1 -1,i_2,\dots,i_n)}(q)$, or
one term each from
$qg_{(i_1-1,i_2+1,i_3,\dots,i_n)}(q)$ and
$g_{(i_1 -1,i_2,\dots,i_n)}(q)$ and two terms from
$g_{(0,i_1+i_2-1,i_3,\dots,i_n)}(q)$.

Let us denote the relevant compositions by
$I:=(i_1,\dots,i_n)$, $I':=(i_1-1,i_2+1,i_3,\dots,i_n)$,
$I'':=(i_1-1,i_2,\dots,i_n)$ and $I''':=(0,i_1+i_2-2,i_3,\dots,i_n)$.

First, consider the terms of $g_{(i_1,\dots,i_n)}(q)$ corresponding to
compositions $J$ such that the first part of $J$ is $i_1$, \emph{i.e.}, $J$
has the form $(i_1,j_2,j_3,\dots,j_r)$.
Let us compare this term to the terms of
$q g_{(i_1-1,i_2+1,i_3,\dots,i_n)}(q)$ and
$g_{(i_1-1,i_2,\dots,i_n)}(q)$ corresponding
to the partitions $J':=(i_1-1,j_2+1,j_3,\dots,j_r)$ and
$J'':=(i_1-1,j_2,j_3,\dots,j_r)$, respectively.
All three terms have the same sign and the same $\stb$:
$\stb(I,J)=\stb(I',J')=\stb(I'',J'')$.
And now it is easy to see that
$q \QStat_B(J')+\QStat_B(J'') = \QStat_B(J)$:
\begin{equation}
\begin{split}
&
q [r]^{i_1-1} [r-1]^{j_2+1} [r-2]^{j_3} \dots 
+ [r]^{i_1-1} [r-1]^{j_2} [r-2]^{j_3} \dots\\
&=
(q [r-1] +1)([r]^{i_1-1} [r-1]^{j_2} [r-2]^{j_3} \dots)\\
& =
  [r]^{i_1} [r-1]^{j_2} [r-2]^{j_3} \dots.
\end{split}
\end{equation}
Note that all terms contain the extra factor $\prod_{t=1}^{r-1} (1+q^t)$.
Therefore the term corresponding to $J$ is equal to the sum of the terms
corresponding to $J'$ and $J''$.

Now consider each term of $g_{(i_1,\dots,i_n)}(q)$ which corresponds to a
composition $J$ such that the first part of $J$ is {\it not} $i_1$,
\emph{i.e.}, $J$ has the form $(j_1,j_2,j_3,\dots,j_r)$ where
$j_1 = i_1+i_2+\dots +i_k$ where $k \geq 2$.
Let us compare this to the following four terms:
the term of $q g_{(i_1-1,i_2+1,i_3,\dots,i_n}(q)$ corresponding
to the composition $J':=J$;
the term of  $g_{(i_1 -1,i_2,\dots,i_n)}(q)$ corresponding
to the composition $J'':=(j_1-1,j_2,\dots,j_r)$;
and the two terms of $g_{(0,i_1+i_2-1,i_3,\dots,i_n)}(q)$
corresponding to the compositions
$J''':=(0,j_1-1,j_2,\dots,j_r)$ and $J^{(4)}:=J''$.
Note that the terms corresponding to
$J, J', J''$, and $J^{(4)}$  have the same sign, while
the term corresponding to $J'''$ has the opposite sign.
And all five terms have the same $\stb$ statistic.
The quantity $\QStat_B$ is nearly the same for every term, and
if we divide each term by $\QStat_B(J'')$, it remains to
verify the equation:
$[r+1]_q = q[r+1]_q+1-(1+q^{r+1}) +1$.
This is clearly true.

We have now accounted for all terms involved in the recurrence.
This completes the proof of the theorem.
\qed

It is very likely that colored Hopf algebra analogues of $\WQSym$, $\FQSym$,
$\NCSF$ already defined in~\cite{MR,Poi,NT,BH,NTqthooks,NTcolored} could be
used to justify the $q$-enumeration of permutation tableaux of type $B$.
Based on preliminary calculations, 
we believe that the type B analogue of 
the matrices from Section \ref{FirstExamples}
are given by computing the transition matrix
between the $S^I$ and two new bases  $\Psi^B_I$
and $L^B(q)$.  Here 
\begin{equation}
\Psi^B_I := \frac{\Psi_I}{\prod_{i=2}^{r}(1+q^{i-1})}
\end{equation}
and $L^B(q)$ is defined by having $M_{L(q),\Psi}$
as transition matrix from the $\Psi^B$.
Note also that this interpretation would
immediately generalize to colored algebras with any number of colors and not
only to two colors.

%%%%%%%%%%%%%%%%%%%%%%%%%%%%%%%%%%%%%%%%%%%%%%%%%%%%%%%%%%%%%%%%%%%%%%%%%%%%%%%
%%%%%%%%%%%%%%%%%%%%%%%%%%%%%%%%%%%%%%%%%%%%%%%%%%%%%%%%%%%%%%%%%%%%%%%%%%%%%%%
%%%%%%%%%%%%%%%%%%%%%%%%%%%%%%%%%%%%%%%%%%%%%%%%%%%%%%%%%%%%%%%%%%%%%%%%%%%%%%%
\section{Appendix -- Conjectures}

We define the {\it descent tops} (also called the 
{\it Genocchi descent set}) of a permutation $\sigma\in S_n$
as $\GDes(\sigma):=\{i\in [2,n] \ \vert \ \sigma(j)=i \Rightarrow
\sigma(j+1) < \sigma(j) \}$.  In other words, 
$\GDes(\sigma)$ is the set of values of the descents of $\sigma$.
We also define the {\it Genocchi composition of descents}
$\GC(\sigma)$ as the integer composition $I$ of $n$ whose descent
set is $\{d-1 \ \vert \ d\in \GDes(\sigma) \}$.

>From Theorem 5.1 of~\cite{HNTT}, it is easy to see that $E_I^J(1)$ is equal to
the number of packed words $w$ such that
\begin{equation}
\GC(\Std(w))=I \qquad\text{and}\qquad \ev(w)=J.
\end{equation}
So $E_I^J(q)$ is the generating function of a statistic in $q$ over this set
of words. We propose the following conjecture:

\begin{conjecture}
\label{conj-initx}
Let $I$ and $J$ be compositions of $n$ and let $W''(I,J)$ be the set of packed
words $w$ such that
\begin{equation}
\GC(\Std(w))=I
\qquad\text{and}\qquad
\ev(w)=J.
\end{equation}
Then
\begin{equation}
E_I^J(q) = \sum_{w\in W''(I,J)} q^{\ttog(w)},
\end{equation}
where $\ttog$ is the  number of occurrences of the patterns $21-1$ and $31-2$
in $w$.
\end{conjecture}

{\footnotesize
For example, the coefficient $2+2q+q^2$ in row $(3,1)$ and column $(2,1,1)$
comes from the fact that the five words
$1132$, $1231$, $1312$, $2311$, and $3112$ respectively have $0$, $0$, $1$,
$1$, and $2$ occurrences of the previous patterns.
}

There should exist a connection between the $\sinv$ statistic and the pattern
counting on special packed words but we have not been able to find it.

Note that packed words $w$ are in bijection with pairs $(\sigma,J)$ where
$\sigma$ is a permutation and $J$ a composition finer than the recoil
composition of $\sigma$. Since the patterns $31-2$ in $\Std(w)$  come from
patterns $21-1$ or $31-2$ in $w$, Conjecture~\ref{conj-initx} is equivalent to

\begin{conjecture}
\label{cor-E}
Let $I$ and $J$ be compositions of $n$ and let $P''(I,J)$ be the set of
permutations $\sigma$ such that
\begin{equation}
\GC(\sigma)=I
\qquad\text{and}\qquad
\DesC(\sigma^{-1})\fatter J.
\end{equation}
Then
\begin{equation}\label{pattern-equation}
E_I^J(q) = \sum_{w\in P'(I,J)} q^{\tto(\sigma)},
\end{equation}
where $\tto(\sigma)$ is the number of occurrences of the
pattern $31-2$ in $\sigma$.
\end{conjecture}

If we apply Sch\"utzenberger's involution to permutations, that is,
$\sigma\mapsto \omega\sigma\omega$, where $\omega=n\cdots 21$
(also known as taking the {\it reverse complement}),
the statistic descent tops is transformed into descent bottoms,
and patterns $31-2$ are transformed into patterns $2-31$.
In that case it follows from results of \cite{SW} that for $J={1^n}$, the sum
in equation (\ref{pattern-equation}) gives the $q$-enumeration of permutation
tableaux of a given shape. 

Therefore if we assume Conjecture~\ref{cor-E},  Theorem~\ref{Theorem2} implies
the following.

\begin{conjecture}
\label{ptai}
When $K=1^n$,
\begin{equation}
E_I^K(q) = \PT^A_I(q) =
\sum_{J\fatter I} (-1/q)^{l(I)-l(J)} q^{-\stb(I,J)} \QStat_A(J).
\end{equation}
\end{conjecture}

Going from $S(q)$ to $R(q)$ is simple, and allows us to reformulate
Conjecture~\ref{ptai} as follows:

\begin{conjecture}
Let $I$ and $J$ be two compositions of $n$.
Let $\PP'(I,J)$ be the set of permutations $\sigma$ such that
$\GC(\sigma)=I$ and $\DesC(\sigma^{-1})=J$.

Then
\begin{equation}
F_I^J(q) = \sum_{\sigma\in \PP(I,J)} q^{\tto(\sigma)}.
\end{equation}
\end{conjecture}

{\footnotesize
For example, the coefficient $1+q+q^2$ in row $(3,1)$ and column $(3,1)$ comes
from the fact that the words $1243$, $1423$, $4123$ respectively have $0$, $1$
and $2$ occurrences of the pattern $31-2$.
}

%%%%%%%%%%%%%%%%%%%%%%%%%%%%%%%%%%%%%%%%%%%%%%%%%%%%%%%%%%%%%%%%%%%%%%%%%%%%%%%
%%%%%%%%%%%%%%%%%%%%%%%%%%%%%%%%%%%%%%%%%%%%%%%%%%%%%%%%%%%%%%%%%%%%%%%%%%%%%%%
%%%%%%%%%%%%%%%%%%%%%%%%%%%%%%%%%%%%%%%%%%%%%%%%%%%%%%%%%%%%%%%%%%%%%%%%%%%%%%%
\section{Tables}
\label{sec-patt}

Here are the transition matrices from $R(q)$ to $L(q)$ (the matrices of the
coefficients $F_I^J(q)$) for $n=3$ and $n=4$,
where the numbers have been replaced by the corresponding list of permutations
having given recoil composition and $\Lig$-composition.

To save space and for better readability, $0$ has been omitted.

\begin{equation}
%{\goth M}''_3 =
%\left(
\begin{array}{|c||c|c|c|c|}
\hline
\text{\rm $\Lig\backslash Rec$}& 3 & 21 & 12 & 111 \\[.1cm]
\hline
\hline
3   & 123 &   &   &   \\[.1cm]
\hline
21  &     & \empd{132}{312} & 213 &   \\[.1cm]
\hline
12  &     &   & 231 &   \\[.1cm]
\hline
111 &     &   &   & 321 \\[.1cm]
\hline
\end{array}
%\right)
\end{equation}

\begin{equation}
%{\goth M}''_4 =
%\left(
\begin{array}{|c||c|c|c|c|c|c|c|c|}
\hline
\text{\rm $\Lig\backslash Rec$}& 4 & 31 & 22 & 211 & 13 & 121 & 112 & 1111
\\[.1cm]
\hline
\hline
4& 1234 &   &   &   &   &   &   &   \\[.1cm]
\hline
31&  & \empd{1243,\ 1423}{4123} & \empd{1324}{3124} &   & 2134 & 2143 &   &  
\\[.1cm]
\hline
22&     &   & \empd{1342}{3142} &   & 2314 &   &   &   \\[.1cm]
\hline
211&     &   & 3412 & \empd{1432,\ 4132}{4312} &   & \empd{2413}{4213} & 3214
&  
\\[.1cm]
\hline
13&     &   &   &   & 2341 &   &   &   \\[.1cm]
\hline
121&     &   &   &   &   & \empd{2431}{4231} & 3241 &   \\[.1cm]
\hline
112&     &   &   &   &   &   & 3421 &   \\
\hline
1111&     &   &   &   &   &   &   & 4321 \\
\hline
\end{array}
%\right)
\end{equation}

%%%%%%%%%%%%%%%%%%%%%%%%%%%%%%%%%%%%%%%%%%%%%%%%%%%%%%%%%%%%%%%%%%%%%%%%%%%%%%%
%%%%%%%%%%%%%%%%%%%%%%%%%%%%%%%%%%%%%%%%%%%%%%%%%%%%%%%%%%%%%%%%%%%%%%%%%%%%%%%
%%%%%%%%%%%%%%%%%%%%%%%%%%%%%%%%%%%%%%%%%%%%%%%%%%%%%%%%%%%%%%%%%%%%%%%%%%%%%%%
%%%%%%%%%%%%%%%%%%%%%  BIBLIOGRAPHIE !!! %%%%%%%%%%%%%%%%%%%%%%%%%%%%%%%%%%%%%%
%%%%%%%%%%%%%%%%%%%%%%%%%%%%%%%%%%%%%%%%%%%%%%%%%%%%%%%%%%%%%%%%%%%%%%%%%%%%%%%
%%%%%%%%%%%%%%%%%%%%%%%%%%%%%%%%%%%%%%%%%%%%%%%%%%%%%%%%%%%%%%%%%%%%%%%%%%%%%%%

\end{document}